\documentclass[12pt, a4paper, oneside]{amsart}
\usepackage[utf8]{inputenc} 
\usepackage{amsmath, amsthm, amssymb, amsfonts, amsthm, mathrsfs} 
\usepackage{cite}
\usepackage{mathtools}
\usepackage{graphicx}
\usepackage{longtable}
\usepackage{eucal}
\usepackage{enumerate}
\usepackage{tikz}             
													
\usepackage{quiver}
\usepackage{booktabs} 		

\usepackage{geometry}					
\usepackage{setspace}
\usepackage{xcolor}

\usepackage{lscape}

\usepackage{memhfixc}

\usepackage[nobysame,alphabetic,msc-links]{amsrefs}

\newcommand{\nc}{\newcommand}
\newcommand{\rc}{\renewcommand}
\nc{\on}{\operatorname}

\nc{\mkcom}[1]{\noindent  \textcolor{red}{$[\star$ M: #1 $\star]$}}
\nc{\gncom}[1]{\noindent  \textcolor{purple}{$[\star$ N: #1 $\star]$}}
\nc{\bwcom}[1]{\noindent  \textcolor{blue}{$[\star$ B: #1 $\star]$}}
\nc{\todo}[1]{\noindent  \textcolor{red}{$[\star$ To do: #1 $\star]$}}

\rc{\AA}{\mathbb{A}}	
\nc{\BB}{\mathbb{B}}	
\nc{\C}{\mathbb{C}}	
\nc{\DD}{\mathbb{D}}	
\nc{\EE}{\mathbb{E}}	
\nc{\F}{\mathbb{F}}	
\nc{\GG}{\mathbb{G}}	
\nc{\HH}{\mathbb{H}}	
\nc{\II}{\mathbb{I}}	
\nc{\JJ}{\mathbb{J}}	
\nc{\KK}{\mathbb{K}}	
\nc{\LL}{\mathbb{L}}	
\nc{\MM}{\mathbb{M}}	
\nc{\NN}{\mathbb{N}}	
\nc{\PP}{\mathbb{P}}	
\nc{\QQ}{\mathbb{Q}}	
\nc{\R}{\mathbb{R}}	
\rc{\SS}{\mathbb{S}}	
\nc{\TT}{\mathbb{T}}	
\nc{\UU}{\mathbb{U}}	
\nc{\VV}{\mathbb{V}}	
\nc{\WW}{\mathbb{W}}	
\nc{\XX}{\mathbb{X}}	
\nc{\YY}{\mathbb{Y}}	
\nc{\Z}{\mathbb{Z}}	

\nc{\bA}{\mathbf{A}}	
\nc{\bB}{\mathbf{B}}	
\nc{\bC}{\mathbf{C}}	
\nc{\bD}{\mathbf{D}}	
\nc{\bE}{\mathbf{E}}	
\nc{\bF}{\mathbf{F}}	
\nc{\bG}{\mathbf{G}}	
\nc{\bH}{\mathbf{H}}	
\nc{\bI}{\mathbf{I}}	
\nc{\bJ}{\mathbf{J}}	
\nc{\bK}{\mathbf{K}}	
\nc{\bL}{\mathbf{L}}	
\nc{\bM}{\mathbf{M}}	
\nc{\bN}{\mathbf{N}}	
\nc{\bO}{\mathbf{O}}	
\nc{\bP}{\mathbf{P}}	
\nc{\bQ}{\mathbf{Q}}	
\nc{\bR}{\mathbf{R}}	
\nc{\bS}{\mathbf{S}}	
\nc{\bT}{\mathbf{T}}	
\nc{\bU}{\mathbf{U}}	
\nc{\bV}{\mathbf{V}}	
\nc{\bW}{\mathbf{W}}	
\nc{\bX}{\mathbf{X}}	
\nc{\bY}{\mathbf{Y}}	
\nc{\bZ}{\mathbf{Z}}	

\nc{\calA}{\mathcal{A}}	
\nc{\calB}{\mathcal{B}}	
\nc{\calC}{\mathcal{C}}	
\nc{\calD}{\mathcal{D}}	
\nc{\calE}{\mathcal{E}}	
\nc{\calF}{\mathcal{F}}	
\nc{\calG}{\mathcal{G}}	
\nc{\calH}{\mathcal{H}}	
\nc{\calI}{\mathcal{I}}	
\nc{\calJ}{\mathcal{J}}	
\nc{\calK}{\mathcal{K}}	
\nc{\calL}{\mathcal{L}}	
\nc{\calM}{\mathcal{M}}	
\nc{\calN}{\mathcal{N}}	
\nc{\calO}{\mathcal{O}}	
\nc{\calP}{\mathcal{P}}	
\nc{\calQ}{\mathcal{Q}}	
\nc{\RR}{\mathcal{R}}	
\nc{\calS}{\mathcal{S}}
\nc{\calT}{\mathcal{T}}	
\nc{\calU}{\mathcal{U}}	
\nc{\calV}{\mathcal{V}}	
\nc{\calW}{\mathcal{W}}
\nc{\calX}{\mathcal{X}}	
\nc{\calY}{\mathcal{Y}}	
\nc{\calZ}{\mathcal{Z}}

\nc{\fraka}{\mathfrak{a}}
\nc{\bb}{\mathfrak{b}}
\nc{\frakc}{\mathfrak{c}}
\nc{\frakd}{\mathfrak{d}}
\nc{\frake}{\mathfrak{e}}
\nc{\frakf}{\mathfrak{f}}
\nc{\frakg}{\mathfrak{g}}
\nc{\frakh}{\mathfrak{h}}
\nc{\fraki}{\mathfrak{i}}
\nc{\frakj}{\mathfrak{j}}
\nc{\frakk}{\mathfrak{k}}
\nc{\frakl}{\mathfrak{l}}
\nc{\frakm}{\mathfrak{m}}
\nc{\frakn}{\mathfrak{n}}
\nc{\frako}{\mathfrak{o}}
\nc{\frakp}{\mathfrak{p}}
\nc{\frakq}{\mathfrak{q}}
\nc{\frakr}{\mathfrak{r}}
\nc{\fraks}{\mathfrak{s}}
\nc{\frakt}{\mathfrak{t}}
\nc{\fraku}{\mathfrak{u}}
\nc{\frakv}{\mathfrak{v}}
\nc{\frakw}{\mathfrak{w}}
\nc{\frakx}{\mathfrak{x}}
\nc{\fraky}{\mathfrak{y}}
\nc{\frakz}{\mathfrak{z}}

\nc{\frakA}{\mathfrak{A}}
\nc{\frakB}{\mathfrak{B}}
\nc{\frakC}{\mathfrak{C}}
\nc{\frakD}{\mathfrak{D}}
\nc{\frakE}{\mathfrak{E}}
\nc{\frakF}{\mathfrak{F}}
\nc{\frakG}{\mathfrak{G}}
\nc{\frakH}{\mathfrak{H}}
\nc{\frakI}{\mathfrak{I}}
\nc{\frakJ}{\mathfrak{J}}
\nc{\frakK}{\mathfrak{K}}
\nc{\frakL}{\mathfrak{L}}
\nc{\frakM}{\mathfrak{M}}
\nc{\frakN}{\mathfrak{N}}
\nc{\frakO}{\mathfrak{O}}
\nc{\frakP}{\mathfrak{P}}
\nc{\frakQ}{\mathfrak{Q}}
\nc{\frakR}{\mathfrak{R}}
\nc{\frakS}{\mathfrak{S}}
\nc{\T}{\mathfrak{t}}
\nc{\frakU}{\mathfrak{U}}
\nc{\frakV}{\mathfrak{V}}
\nc{\frakW}{\mathfrak{W}}
\nc{\frakX}{\mathfrak{X}}
\nc{\frakY}{\mathfrak{Y}}
\nc{\frakZ}{\mathfrak{Z}}

\nc{\sA}{\mathscr{A}}	
\nc{\sB}{\mathscr{B}}	
\nc{\sC}{\mathscr{C}}	
\nc{\sD}{\mathscr{D}}	
\nc{\sE}{\mathscr{E}}	
\nc{\sF}{\mathscr{F}}	
\nc{\sG}{\mathscr{G}}	
\nc{\sH}{\mathscr{H}}	
\nc{\sI}{\mathscr{I}}	
\nc{\sJ}{\mathscr{J}}	
\nc{\sK}{\mathscr{K}}	
\nc{\sL}{\mathscr{L}}	
\nc{\sM}{\mathscr{M}}	
\nc{\sN}{\mathscr{N}}	
\nc{\sO}{\mathscr{O}}	
\nc{\sP}{\mathscr{P}}	
\nc{\sQ}{\mathscr{Q}}	
\nc{\sR}{\mathscr{R}}	
\nc{\sS}{\mathscr{S}}
\nc{\sT}{\mathscr{T}}	
\nc{\sU}{\mathscr{U}}	
\nc{\sV}{\mathscr{V}}	
\nc{\sW}{\mathscr{W}}
\nc{\sX}{\mathscr{X}}	
\nc{\sY}{\mathscr{Y}}	
\nc{\sZ}{\mathscr{Z}}

\nc{\Lie}{\on{Lie}}
\nc{\GL}{\on{GL}}
\nc{\PGL}{\on{PGL}}
\nc{\PSL}{\on{PSL}}
\nc{\SL}{\on{SL}}
\nc{\Sp}{\on{Sp}}
\nc{\GSp}{\on{GSp}}
\nc{\GU}{\on{GU}}
\nc{\SO}{\on{SO}}
\nc{\Or}{\on{O}}
\nc{\gl}{\on{\mathfrak{gl}}}
\rc{\sl}{\on{\mathfrak{sl}}}
\nc{\Mat}{\on{Mat}}
\nc{\Fun}{\on{Fun}}
\nc{\Aut}{\on{Aut}}
\nc{\End}{\on{End}}
\nc{\Hom}{\on{Hom}}
\nc{\Mor}{\on{Mor}}
\nc{\Sym}{\on{Sym}}
\nc{\Span}{\on{Span}}
\nc{\Irr}{\on{Irr}}
\nc{\Uch}{\on{Uch}}
\nc{\Type}{\on{Type}}
\nc{\Spec}{\on{Spec}}
\nc{\Ind}{\on{Ind}}
\nc{\Res}{\on{Res}}
\nc{\stab}{\on{stab}}
\nc{\orb}{\on{orb}}
\rc{\ker}{\on{ker}}
\nc{\im}{\on{im}}
\nc{\tr}{\on{tr}}
\nc{\ord}{\on{ord}}
\nc{\Tor}{\on{Tor}}
\nc{\Ad}{\on{Ad}}
\nc{\Id}{\on{Id}}
\nc{\Log}{\on{Log}}
\nc{\Exp}{\on{Exp}}
\nc{\Frac}{\on{Frac}}
\nc{\diag}{\on{diag}}
\nc{\Deck}{\on{Deck}}
\nc{\D}{\on{D}}

\nc{\Nil}{\on{Nil}}
\nc{\Max}{\on{Max}}
\nc{\rank}{\on{rank}}
\nc{\Spin}{\on{Spin}}

\nc{\St}{\mathrm{St}}
\nc{\triv}{\mathrm{triv}}
\nc{\sgn}{\mathrm{sgn}}
\nc{\reg}{\mathrm{reg}}
\nc{\op}{\mathrm{op}}
\nc{\ad}{\mathrm{ad}}
\nc{\id}{\mathrm{id}}
\rc{\ss}{\mathrm{ss}}
\nc{\HLV}{\mathrm{HLV}}
\nc{\Sch}{\mathrm{Sch}}
\nc{\an}{\mathrm{an}}
\nc{\red}{\mathrm{red}}
\rc{\|}{|\!|}

\DeclareMathOperator{\Stab}{Stab}

\DeclareMathOperator{\SAut}{\text{SAut}}

\theoremstyle{plain} 
\newtheorem{theorem}{Theorem}
\newtheorem{lemma}[theorem]{Lemma}
\newtheorem{proposition}[theorem]{Proposition}

\newtheorem{corollary}[theorem]{Corollary}

\theoremstyle{definition} 
\newtheorem{definition}[theorem]{Definition}
\newtheorem{example}[theorem]{Example}

\theoremstyle{remark} 
\newtheorem{remark}[theorem]{Remark}

\makeatletter
\newcommand{\sectionnotoc}[1]{%
  \begingroup
  \let\@tocwrite\@gobbletwo
  \section*{#1}
  \endgroup
}
\makeatother

\title[Holomorphic flexibility properties of strongly commuting varieties]{Holomorphic flexibility properties\\ of strongly commuting varieties}
\author{Stefano Giannini}
\address{School of Mathematical Sciences, College of Science,
Adelaide University, SA 5005, Australia}
\email{stefano.giannini@adelaide.edu.au}

\begin{document}

\begin{abstract}
We study holomorphic flexibility properties of the strongly commuting variety of a complex semisimple algebraic group $G$. We prove that the smooth locus of the strongly commuting variety of $G$ is Oka-$1$. We also prove that the smooth locus of the commuting variety of a semisimple complex Lie algebra is Oka-$1$, and that it contains a Zariski dense open subset which is algebraically elliptic.
\end{abstract}

\maketitle

\tableofcontents

\section{Introduction}

\makeatletter
\begingroup
\let\@tocwrite\@gobbletwo

\subsection*{Main results}
This paper studies holomorphic flexibility properties of two commuting varieties: the commuting variety of a complex algebraic group $G$, defined as
\[R(G) := \{(x,y)\in G^2\,|\,xy=yx\},\]
and the commuting variety of a complex Lie algebra $\frakg$, defined as
\[\calC(\frakg):=\{(x,y)\in \frakg^2\,|\, [x,y]=0\}.\]

The Oka-$1$ property is a holomorphic flexibility property introduced by Alarc\'on and Forstneri\v c \cite{AF25}. Roughly speaking, a complex manifold is Oka-$1$ if every continuous map from an open Riemann surface into it is homotopic to a holomorphic map, with approximation and interpolation. The precise definition is recalled in \S~\ref{c:prelimsOka}.

For a semisimple algebraic group $G$, let $q:\widetilde{G}\to G$ be its universal covering. Define 
\[R_{\rm str}(G) :=(q\times q)(R(\widetilde{G})) \]
and call it the strongly commuting variety of $G$. Thus, a pair belongs to $R_{\rm str}(G)$ precisely when its entries lift to commuting elements in $\widetilde{G}$. Observe that if $G$ is simply connected then $R(G) = R_{\rm str}(G)$. Our main results are the following:
\begin{theorem}\label{thm:theorem1}
Let $G$ be a semisimple complex algebraic group. The smooth locus of the strongly commuting variety $R_{\rm str}(G)$ of $G$ is Oka-$1$.
\end{theorem}
In particular, Theorem~\ref{thm:theorem1} implies that any continuous family of commuting pairs in the smooth locus of $R_{\rm str}(G)$, parameterised by an open Riemann surface, is homotopic to a holomorphic family, with approximation and interpolation. We prove an analogous result in the Lie algebra case:
\begin{theorem}\label{thm:theorem2}
Let $\frakg$ be a semisimple complex Lie algebra. The smooth locus of the commuting variety $\calC(\frakg)$ of $\frakg$ is Oka-$1$. Moreover, there is a dense open subset of $\calC(\frakg)$ which is algebraically elliptic.
\end{theorem}
\par
\endgroup
\makeatother

\makeatletter
\begingroup
\let\@tocwrite\@gobbletwo
\subsection*{Historical background}

Commuting varieties were first studied in \cite{MT55} for the Lie algebra $\gl_n$. The authors showed that the commuting variety is irreducible of dimension $n^2+n$. This result was later generalised by Richardson \cite{Ric79} who proved that the commuting varieties of reductive Lie algebras (resp. semisimple simply connected groups) over algebraically closed fields of characteristic zero are irreducible. The Lie algebra result was extended in \cite{Lev02} in the case of a field of positive characteristic and under certain mild assumptions on $G$. Other versions of the commuting variety have also been studied. For instance, in \cite{Bar01} and \cite{Pre03}, motivated by the study of moduli of vector bundles on algebraic surfaces, the authors investigated commuting varieties of nilpotent endomorphisms of a vector space, and related them to certain punctual Hilbert schemes. For other variants, in \cite{BFM02} and \cite{KS00} the authors studied and classified varieties of almost commuting pairs of elements in compact Lie groups $K$ which appear naturally when one studies flat connections on principal $K$-bundles over a torus. One may also consider the variety of commuting $n$-tuples of a reductive group $G$ and of its maximal compact subgroup $K$. The topology and cohomology of these varieties have been studied in \cite{Bai07,GPS12,PS13}. In particular, it was shown that the variety of commuting $n$-tuples in $K$ is a strong deformation retract of the variety of commuting $n$-tuples in $G$. Commuting varieties of algebraic groups are also special cases of representation varieties.

\medskip

Representation varieties $\Hom(\pi_1(\Sigma),G)$ are varieties of representations of the fundamental group of a surface $\Sigma$ into a group $G$. Thus, the commuting variety of $G$ is the representation variety associated with a torus and the group $G$. The character variety associated with $\Sigma$ and $G$ is the quotient of the representation variety by the conjugation action of $G$, parameterising representations up to equivalence. Character varieties, together with their related moduli spaces of flat connections and Higgs bundles, have been studied using a wide range of methods, appearing in many areas of mathematics and physics including non-abelian Hodge theory \cite{Sim92,Sim94}, mirror symmetry \cite{HT03}, geometric Langlands \cite{BD96, BZN18}, and gauge theory \cite{AB83,Hit87}. In this paper we study commuting varieties from the point of view of Oka theory.
\medskip

Hyperbolic manifolds and Oka manifolds sit at opposite ends in the landscape of complex manifolds. The former are holomorphically \emph{rigid} and the latter are holomorphically \emph{flexible}. For instance, a theorem of Liouville implies that holomorphic maps $\C\to \mathbb D$ are constant, illustrating the rigidity of such maps. Oka theory makes precise the property of a complex manifold of admitting `many' holomorphic maps into it. Oka theory originates in the classical work of Oka and Grauert concerning analytic problems on Stein manifolds. This work gave rise to the Oka principle which says that analytic problems on Stein spaces have analytic solutions whenever there are no topological obstructions. In his seminal 1989 paper \cite{Gro89}, Gromov introduced the notion of a dominating spray, a geometric condition which implies strong holomorphic flexibility properties. More precisely, he proved that if $X$ admits a dominating spray, then $X$ satisfies an Oka principle for holomorphic maps from Stein spaces $S$: every continuous map $S\to X$ is homotopic to a holomorphic one with approximation and interpolation.  Building on Gromov's work, Forstneri\v c and others developed the modern theory of Oka manifolds and established the equivalence of several approximation, interpolation, and parametric forms of the Oka property \cite{For06,For17}.

Oka theory has now seen a diverse range of applications. These include the classification of holomorphic vector bundles on Stein spaces, the existence of holomorphic sections with prescribed local properties, problems concerning holomorphic embeddings, submersions, and matrix factorisations \cite{For18, IK12}. An application closer to the setting of this paper appears in the work of Guralnick \cite{Gur88} and Leiterer \cite{Lei19} where the authors prove that holomorphic maps $A,B:S\to \Mat_n(\C)$ from $1$-dimensional Stein spaces $S$ which are locally holomorphically similar are globally holomorphically similar. Guralnick proved this when $S$ is an open Riemann surface and Leiterer extended it to arbitrary 1-dimensional Stein spaces. In \cite{Lei19}, Leiterer exploited an Oka principle for Oka pairs (in the sense of Forster and Ramspott \cite{FR66}) to reduce the assertion that locally holomorphically similar matrices are globally holomorphically similar to a certain topological statement and then proved the topological statement. In \cite[Theorem 6.14.9]{For17}, the results of Guralnick and Leiterer are obtained via an application of an Oka principle for sections of branched holomorphic maps (cf. \cite[Theorem 6.14.6]{For17}).  

In \cite{AF25}, Alarc{\'o}n and Forstneri\v c introduced the Oka-$1$ property which restricts the source manifolds to open Riemann surfaces. Thus, every Oka manifold is Oka-$1$. Roughly speaking, Oka-$1$ manifolds are characterised by the property that they admit many holomorphic maps from open Riemann surfaces. They proved that a complex manifold dominable by $\C^n$ at most points is Oka-$1$. Furthermore, they introduced a local spray approximation property which implies the Oka-$1$ property, and is one of the main tools used in this paper.

\par
\endgroup
\makeatother

\makeatletter
\begingroup
\let\@tocwrite\@gobbletwo
\subsection*{Idea of proof}
We outline the key ideas of the proof of Theorem~\ref{thm:theorem1}; the proof of Theorem~\ref{thm:theorem2} is analogous. Suppose first that $G$ is simply connected, and put $r=\rank(G)$. The commuting variety $R(G)$ is singular, so flexibility cannot be applied directly to the whole variety. Instead, we consider the open subset of pairs whose first entry is regular, and show that it lies in the smooth locus. The fibre of the first projection map to $G$ over a regular element $x$ of $G$ is its centraliser $Z_G(x)$. Although the centraliser is abelian, it is not generally connected, so the exponential map from its Lie algebra only hits the identity component. We prove a structural result for $Z_G(x)$ (Lemma~\ref{lem:regular-centraliser-components}) which says the remaining components can be reached by multiplication by powers of $x$ and central elements. This leads us to a family of local biholomorphisms 
\[\Phi_{(n,\varepsilon)}:G_{\rm reg}\times\C^r
\longrightarrow R_{\rm reg}(G),
\qquad (n,\varepsilon)\in\Z\times Z(G),\]
whose images cover $R_{\rm reg}(G)$. These maps are constructed using differentials of the fundamental characters of $G$ \cite[\S 8]{Ste65}. The source $G_{\rm reg}\times \C^r$ is flexible and hence has the local spray approximation property (LSAP), recalled in \S\ref{c:prelimsOka}. The local biholomorphisms transfer this property to $R_{\rm reg}(G)$. The complement of $R_{\rm reg}(G)$ in the smooth locus of $R(G)$ is an analytic subset of positive codimension, so a criterion \cite[Proposition 7.13(b)]{AF25} implies the whole smooth locus is Oka-$1$.

For a general semisimple group $G$, we pass to the universal covering $q:\widetilde{G}\to G$. Since $R(\widetilde{G})$ is the inverse image of $R_{\rm str}(G)$ under $q\times q$, restricting this covering to $R_{\rm reg}(\widetilde{G})$ allows us to transfer LSAP from $R_{\rm reg}(\widetilde{G})$ to the regular strongly commuting variety of $G$. Again, using the criterion in \cite[Proposition 7.13(b)]{AF25} gives us the Oka-$1$ property for $R_{\rm str}(G)_{\rm sm}$.
\par
\endgroup
\makeatother

\makeatletter
\begingroup
\let\@tocwrite\@gobbletwo

\subsection*{The full Oka property}
Results in the very recent preprint \cite{ADG26}, which originated in \cite{Ope26}, show that our proofs of Theorem~\ref{thm:theorem1} and Theorem~\ref{thm:theorem2} give a much stronger conclusion, namely the full Oka property. Indeed, the proof of Theorem~\ref{thm:theorem1} shows that we can find a family of local biholomorphisms $\Phi_{(n,\varepsilon)}: G_{\rm reg}\times \C^r\to R_{\rm reg}$ whose images cover $R_{\rm reg}(G)$ and where $G_{\rm reg}\times \C^r$ is flexible and hence strongly dominable. It follows that $R_{\rm reg}(G)$ is strongly dominable and hence Oka by \cite[Theorem 1.4]{ADG26}. Moreover, by \cite[Theorem 1.8]{ADG26}, the Oka property extends across analytic subsets of positive codimension. Similarly, for the strongly commuting variety $R_{\rm str}(G)$ of a semisimple group $G$, there is a surjective local biholomorphism $R_{\rm reg}(\widetilde{G})\to R_{\rm str, reg}(G)$ which makes $R_{\rm str,reg}(G)$ strongly dominable and hence Oka by \cite[Theorem 1.4]{ADG26}. Thus, the whole smooth locus $R_{\rm str}(G)_{\rm sm}$ is Oka by \cite[Theorem 1.8]{ADG26}. Similarly, the conclusions that the smooth loci of the commuting variety of a complex semisimple Lie algebra $\frakg$ and the commuting variety of $\PSL_p(\C)$ for prime $p$ are Oka-$1$ can be upgraded to the full Oka property by the same reasoning.
\par
\endgroup
\makeatother

\makeatletter
\begingroup
\let\@tocwrite\@gobbletwo
\subsection*{Outline of the paper}
In \S~\ref{c:prelims} we briefly discuss some preliminaries regarding the Oka-$1$ property and the local spray approximation property which is a key tool in the proofs of the main results. We briefly recall the definition of regular elements in a simply connected semisimple group $G$ and prove a structural result for the centralisers of such elements. Section~\ref{c:regularcommutingvariety} is devoted to the study of the regular commuting variety, a Zariski open subset of the commuting variety of $G$ that plays a central role in the proof. In \S\ref{c:Smoothness}, we prove that the regular commuting variety is contained in the smooth locus of the commuting variety. We then define the local biholomorphisms in \S~\ref{c:localparams} and use these to prove Theorem~\ref{thm:theorem1} in \S~\ref{c:Oka1theorem} for the case of simply connected semisimple groups. In \S~\ref{c:stronglycommutingvarieties}, we define the strongly commuting variety and give four equivalent descriptions of which two are topological and two are algebraic. In \S~\ref{c:Oka1theoremstronglycommutingvarieties}, we prove Theorem~\ref{thm:theorem1}. In \S~\ref{c:Oka1PSLp}, we drop the simple connectedness assumption by proving the Oka-$1$ property for the smooth locus of the commuting variety of $\PSL_p(\C)$ when $p$ is prime. Finally, in \S~\ref{c:Liealgebra}, we move to the Lie algebra case and prove Theorem~\ref{thm:theorem2}.
\par
\endgroup
\makeatother

\makeatletter
\begingroup
\let\@tocwrite\@gobbletwo
\subsection*{Acknowledgements}
The author would like to thank Finnur L\'arusson for the many helpful discussions, for carefully reading the manuscript and providing feedback, and for pointing out the upgrade from the Oka-$1$ property to the full Oka property. The author is also grateful to David Baraglia for reading and providing helpful comments on this paper. This research was supported by the Australian Government Research Training Program.
\par
\endgroup
\makeatother

\makeatletter
\begingroup
\let\@tocwrite\@gobbletwo
\subsection*{Use of artificial intelligence}
This paper was written entirely by the author. ChatGPT was used for finding literature, and GPT-5.6 was used during the preparation of the manuscript for proof-reading and checking mathematical arguments. In particular, it helped identify minor errors in Lemma~\ref{lem:regular-centraliser-components} and Theorem~\ref{thm:RregLSAP}. 
\par
\endgroup
\makeatother

\makeatletter
\begingroup
\let\@tocwrite\@gobbletwo
\subsection*{Conventions}
Throughout this paper, all algebraic groups, varieties and Lie algebras are defined over the field of complex numbers $\C$. We use the convention that semisimple algebraic groups are connected.
\par
\endgroup
\makeatother

\section{Preliminaries}\label{c:prelims}
In this section, we recall the notions of holomorphic flexibility used in the paper and establish some of the key ingredients used in the proof of Theorem~\ref{thm:theorem1}. We give a structural result for the centralisers of regular elements in a simply connected almost simple algebraic group, and recall how differentials of the fundamental characters of $G$ can provide a basis for their Lie algebras.

\subsection{Holomorphic flexibility and Oka-$1$ manifolds}\label{c:prelimsOka}
There are several precise notions of holomorphic flexibility. Here, we define the notions relevant to the paper. These notions pertain to smooth varieties and they include flexibility, (algebraic) ellipticity, the local spray approximation property, and the Oka-$1$ property. 

The notion of a flexible variety was introduced by Arzhantsev et al. \cite{AZK12} and further developed in \cite{AFK13}.
\begin{definition}
 An irreducible variety $X$ is \emph{flexible} if for all smooth points $x\in X$, the tangent space $T_xX$ is spanned by the tangent vectors to the orbits $Hx$ of one-parameter unipotent subgroups $H\subset \Aut(X)$.
\end{definition}
For an irreducible affine variety $X$ of dimension at least two, flexibility is equivalent to the subgroup $\SAut(X)\subset \Aut(X)$, generated by one-parameter unipotent subgroups, acting infinitely transitively on the smooth locus of $X$ \cite{AFK13}. Furthermore, it is shown in \cite[Proposition A.3]{AFK13} that if a smooth irreducible variety is flexible then it admits a dominating spray, a geometric condition defined by Gromov \cite{Gro89} which implies the Oka principle for holomorphic sections of holomorphic fibre bundles over Stein manifolds.
\begin{definition}
Let $X$ be a complex manifold. A \emph{dominating spray} on $X$ is a holomorphic
vector bundle $p:E\to X$ together with a holomorphic map $s:E\to X$ such
that $s$ restricts to the identity on the zero section $Z$, and for each $x\in Z\cong X$ the tangent map $d_xs$ sends the fibre $E_x=p^{-1}(x)$ surjectively onto $T_xX$. Furthermore, $X$ is called \emph{elliptic} if it admits a dominating spray. We say that $X$ is \emph{algebraically elliptic} if $X$ is a smooth algebraic variety and the holomorphic vector bundle, together with the dominating spray, are algebraic.
\end{definition}
Every elliptic manifold is Oka (cf. \cite[Corollary 5.6.14]{For17}) and hence Oka-$1$. The Oka-$1$ property defines a class of complex manifolds characterised by the property that they admit many holomorphic curves parameterised by open Riemann surfaces. This class of complex manifolds was introduced in \cite{AF25} and further studied in \cite{FL25}. Let us recall the definition.

\begin{definition}
A connected complex manifold $X$ is an \emph{Oka-$1$ manifold} if for any open Riemann surface $\Sigma$, compact Runge set $K$ in $\Sigma$, a discrete sequence $(a_i)_{i\in I}$ in $\Sigma$ without repetitions, continuous map $f:\Sigma\to X$ which is holomorphic on a neighbourhood of $K\cup \bigcup_i \{a_i\}$, distance function $d$ inducing the topology of $X$, number $\epsilon >0$, and integers $k_i\geq 1$ there is a holomorphic map $F:\Sigma \to X$ which is homotopic to $f$ and satisfies
\begin{enumerate}
\item[(1)] $\sup_{p\in K}d(F(p),f(p))< \epsilon$, and
\item[(2)] $F$ agrees with $f$ to order $k_i$ at the point $a_i$ for every $i$.
\end{enumerate} 
\end{definition}

We now introduce the local spray approximation property (LSAP). This property, introduced in \cite{AF25}, implies the Oka-$1$ condition \cite{AF25}. The LSAP property plays an important role when we study the holomorphic flexibility properties of commuting varieties.

A pair of compact topological discs $D\subset D'\subset \C$ is called a \emph{special pair} if both have piecewise smooth boundaries and $D'\setminus \mathring{D} $ is a disc attached to $D$ along a boundary arc. A holomorphic map $F:D\times \BB^N \to X$ to a complex manifold $X$ is called a holomorphic spray of maps $D\to X$ with core $f= F(\cdot,0)$, where $\BB^N$ denotes the unit ball in $\C^N$. Such a spray is said to be \emph{dominating} if the partial derivative $\frac{\partial}{\partial t}|_{t=0}F(z,t):\C^N\to T_{f(z)}X$ is surjective for all $z\in D$.

\begin{definition}
A complex manifold $X$ has the \emph{local spray approximation property} (LSAP) at $x\in X$ if there is an open neighbourhood $V\subset X$ of $x$ satisfying the following condition. Given a special pair of compact discs $D\subset D'\subset \C$ and dominating holomorphic spray $F:D\times \BB^N\to V$, there is a number $r\in (0,1)$ such that  $F$ can be approximated uniformly on $D\times r \overline{\BB}^N$ by holomorphic maps $G:D'\times  r\overline{\BB}^N \to X$. A manifold $X$ has LSAP if the above condition holds at every $x\in X$. 
\end{definition}

By \cite[Proposition 7.13 (b)]{AF25}, every manifold with LSAP is Oka-$1$. It is not known whether every Oka-$1$ manifold has LSAP.

\subsection{Regular elements and centralisers}

Let $G$ be a semisimple simply connected linear algebraic group. Following \cite{Ste65}, we say an element $x\in G$ is \emph{regular} if its centraliser $Z_G(x)$ in $G$ has minimal dimension equal to $r:=\rank(G)$. In \cite{Ste65}, it is shown that the subset $G_{\mathrm{reg}}$ of regular elements is Zariski open and dense in $G$. Moreover, since $G$ is assumed to be simply connected, $Z_G(x)$ is abelian \cite[Theorem D]{Kur78}.

\subsubsection{Centralisers of regular elements}
Here we give a decomposition of the centraliser in $G$ of a regular element in a simply connected almost simple algebraic group $G$, which is vital for the proof of Theorem~\ref{thm:theorem1}. The corresponding result for simply connected semisimple groups follows by applying the lemma to each almost simple factor. The proof is a consequence of \cite[Remark 3]{Som98}; see also \cite[Lemma 34]{MS03}.
For an algebraic group $H$, we denote by $H^\circ$ the connected component of the identity of $H$.

\begin{lemma}\label{lem:regular-centraliser-components}
Let $G$ be a simply connected, almost simple algebraic group and let $x\in G_{\mathrm{reg}}$ be a regular element. Then
\[Z_G(x)=\bigcup_{ n\in\mathbb Z}x^{n}Z_G(x)^\circ Z(G).\]
Equivalently, the quotient
\[Z_G(x)/(Z_G(x)^\circ Z(G))\]
is finite cyclic generated by the class of $x$.
\end{lemma}
\begin{proof}
Write $x=x_sx_u=x_ux_s$ for the Jordan decomposition of $x$, and let $H:=Z_G(x_s)$ be the centraliser of $x_s$ in $G$.
Using the uniqueness of the Jordan decomposition, it is readily checked that
\[Z_G(x)=Z_G(x_s)\cap Z_G(x_u)=Z_H(x_u).\]
Further, since $x\in G_{\rm reg}$, the unipotent part $x_u$ is a regular unipotent element of $H$ \cite[Section 1.14]{Car93}, and since $G$ is simply connected, \cite[Theorem 2.11]{Hum95} shows that $H$ is a connected reductive subgroup of maximal rank. Since $x_u\in H$ is regular, it lies in a unique Borel subgroup $B=TU$ of $H$, where $T$ is the maximal torus and $U$ the unipotent radical of $B$. It follows from \cite[Proposition 4.1]{Hum95} that $Z_H(x_u) = Z_B(x_u)$ and $Z_B(x_u) = Z(H)Z_U(x_u)$. Thus, 
\[Z_G(x)=Z_H(x_u)=Z(H)Z_U(x_u).\]
Consider the quotient map $\pi:G\to G_{\rm ad}:=G/Z(G)$ and put $\overline{H}=Z_{G_{\rm ad}}(\bar{x}_s)^\circ$. Since $H = Z_G(x_s)$ is connected and $\pi$ is a central isogeny, we have $\pi(H) =\overline{H}$.  Moreover, 
\[\pi_Z:=\pi|_{Z(H)}:Z(H)\longrightarrow Z(\overline{H})\]
is surjective. Indeed, if $h\in H$ and $\pi(h)\in Z(\overline{H})$ then the homomorphism
\[H\longrightarrow Z(G),\quad a\longmapsto [h,a]\]
is constant, because $H$ is connected and $Z(G)$ is finite, and its value at the identity is $1$. Hence $[h,a]=1$ for every $a\in H$, so $h\in Z(H)$.
By \cite[Remark 3]{Som98}, the quotient
\[Z(\overline{H})/Z(\overline{H})^\circ\]
is a finite cyclic group generated by the class $\bar{x}_s\cdot Z(\overline{H})^\circ$. We claim that
\[Z(H)/Z(H)^\circ Z(G) \cong Z(\overline{H})/Z(\overline{H})^\circ.\]
Indeed, consider the map
\[\begin{tikzcd}[ampersand replacement=\&]
	{Z(H)} \& {Z(\overline{H})} \& {Z(\overline{H})/Z(\overline{H})^\circ.}
	\arrow["{\pi_Z}", from=1-1, to=1-2]
	\arrow[from=1-2, to=1-3]
\end{tikzcd}\]
It is surjective with kernel given by $\pi^{-1}_Z( Z(\overline{H})^\circ)=Z(H)^\circ Z(G)$.
Thus,
\[Z(H)=\bigcup_{n\in\Z}x_{s}^{n} \cdot Z(H)^\circ Z(G).\]
Substituting this decomposition into $Z_G(x)=Z(H)Z_U(x_u)$, we obtain
\[\begin{aligned}
        Z_G(x)
        &=Z(H)Z_U(x_u) \\
        &= \bigcup_{n\in\mathbb Z} x_s^n Z(H)^\circ Z_{U}(x_u)Z(G).
\end{aligned}\]
Since $Z(H)^\circ$ and $Z_{U}(x_u)$ are connected subgroups of $Z_{G}(x)^\circ$, their product is contained in $Z_G(x)^\circ$. Therefore
\[Z(H)^\circ Z_{U}(x_u) Z(G)\subset Z_{G}(x)^\circ Z(G).\]
It follows that
\[Z_G(x)\subset \bigcup_{ n\in\Z}x_s^{n} Z_{G}(x)^\circ Z(G).\]
The reverse inclusion is clear. We conclude that
\[Z_G(x)=\bigcup_{ n\in\Z}x_s^{n}Z_G(x)^\circ Z(G).\]
Finally, we may replace $x_s$ by $x$ in this decomposition. Indeed, since $x\,x_s^{-1}=x_u$ and
\[x_{u}\in Z_{U}(x_{u})\subset Z_{G}(x)^\circ,\]
the images of $x$ and $x_{s}$ in the quotient
\[Z_{G}(x)/(Z_{G}(x)^\circ Z(G))\]
coincide. Hence
\[Z_G(x)=\bigcup_{ n\in\Z}x^n Z_G(x)^\circ Z(G).\]
Equivalently, the quotient
\[Z_G(x)/Z_G(x)^\circ Z(G)\]
is generated by the class of $x$.
\end{proof}
\begin{remark}
When $G=\SL_n(\C)$ and $x\in G$ is regular, the quotient 
\[Z_G(x)/Z_G(x)^\circ Z(G)\]
is trivial. Indeed, in the notation of the proof of the lemma, we saw above that 
\[Z_G(x) = Z(H) Z_U(x_u).\]
As $G$ is of type $A$, the pseudo-Levi subgroup $H = Z_G(x_s)$ is actually a Levi subgroup. By \cite[Proposition 3.5]{CE12} this is equivalent to $Z(H)=Z(H)^\circ Z(G)$. We saw above that 
\[Z_G(x) = Z(H)Z_U(x_u) = Z(H)^\circ Z(G)Z_U(x_u)\]
and that $Z(H)^\circ Z_U(x_u)\subset Z_G(x)^\circ$, hence $Z_G(x)\subset Z_G(x)^\circ Z(G)$ and the quotient \[Z_G(x)/Z_G(x)^\circ Z(G)\]
is trivial.
\end{remark}
We provide an example in which this quotient is non-trivial.

\begin{example}
Let $G = \Sp_4(\C) = \{g\in \SL_4(\C)\,|\, g^tAg=A\}$, where 
\[A = \begin{pmatrix}
0&1&0&0\\ -1&0&0&0\\ 0&0&0&1\\0&0&-1&0
\end{pmatrix}.\] 
Consider the semisimple element
\[s = \begin{pmatrix}
1&0&0&0\\ 0&1&0&0\\ 0&0&-1&0\\0&0&0&-1
\end{pmatrix}\in T = \left\{\begin{pmatrix}
t_1&0&0&0\\ 0&t_1^{-1}&0&0\\ 0&0&t_2&0\\0&0&0&t_2^{-1}
\end{pmatrix}\,:t_1,\,t_2\in \C^\times\right\}.\]
The centraliser of $s$ is given by $H=Z_G(s) \simeq \SL_2(\C)\times \SL_2(\C).$
Thus, the centre of $H$ is
\[Z(H) \simeq \mu_2\times \mu_2,\]
where $\mu_n$ denotes the group of $n^{\rm th}$-roots of unity. Thus, the identity component $Z(H)^\circ$ is trivial. However, the centre of $G$ is
\[Z(G) = \{I_2\oplus I_2, -I_2\oplus -I_2\}\simeq \mu_2,\]
where $I_2$ denotes the 2-by-2 identity matrix. Thus,
\[
Z(H)/(Z(H)^\circ Z(G)) =Z(H)/Z(G) \simeq\mu_2,
\]
and the class of $s$ generates this quotient. Choose the regular unipotent element 
\[u =\begin{pmatrix}
1&1&0&0\\ 0&1&0&0\\ 0&0&1&1\\0&0&0&1
\end{pmatrix}\]
of $H$ and let $a=su$. If $U$ is the unipotent radical of the unique Borel subgroup of $H$ containing $u$, then
\[Z_G(a) = Z(H)Z_U(u)\]
and
\[Z_U(u) = \left\{\begin{pmatrix}1&c\\0&1\end{pmatrix}\oplus\begin{pmatrix}1&d\\0&1\end{pmatrix}\in U\bigg|\,c,d\in \C\right\}.\]
In particular, $a$ is regular in $G$ since $\dim (Z_G(a))=\dim (Z_U(u))=2$. Moreover, $Z_G(a)$ has four connected components corresponding to the elements in the centre of $H$,
\[Z(H) = \{(I_2\oplus I_2), (I_2\oplus -I_2),(-I_2\oplus I_2),(-I_2\oplus -I_2)\}.\]
To see that $Z_G(x)\neq Z_G(x)^\circ Z(G)$ in general, let $I_2\oplus I_2$ and $-I_2\oplus -I_2$ be the central elements of $G$. Then the union
\[(I_2\oplus I_2)Z_U(u)\cup (-I_2\oplus -I_2)Z_U(u) = Z_G(a)^\circ Z(G)\]
does not cover all of $Z_G(a)$, as it does not contain the connected components represented by $(I_2\oplus -I_2)$ and $ (-I_2\oplus I_2).$ Since $a=su$ and $u\in Z_G(a)^\circ$, the classes of $a$ and $s$ in
\[Z_G(a)/(Z_G(a)^\circ Z(G))\]
coincide. Thus, multiplication by $a$ gives the remaining coset 
\[ Z_G(a)=Z_G(a)^\circ Z(G)\cup a\cdot Z_G(a)^\circ Z(G).\]
\end{example}

\subsubsection{Differentials of invariant functions}\label{s:differentialsofinvariantfunctions}
Recall that the action of $G$ on itself by conjugation induces an action of $G$ on its ring of functions. We denote by $\C[G]^G$ the corresponding ring of invariants. For $f\in \C[G]^G$, the differential $df_g:T_g G\to \C$, after a translation, defines an element in the dual space $\frakg^*$ given by
\[df_g(X) = \frac{d}{dt}\bigg|_{t=0}f(g\exp(tX)).\]
Since we assume $G$ is semisimple, the Killing form $\kappa$ gives an isomorphism $\frakg \to \frakg^*$. Let $F_g:=\kappa^{-1}(df_g)\in \frakg$ be the element in $\frakg$ corresponding to $df_g$ under this isomorphism. The element $F_g$ is also defined by the formula
\[\kappa(F_g, X) = \frac{d}{dt}\bigg|_{t=0}f(g\exp(tX)),\quad \text{for all}\,\,X\in \frakg.\]
The following is stated in \cite[page 264]{Gol86}; we prove it for completeness:
\begin{proposition}\label{prop:differentialinvariantfunction}
Let $G$ be a semisimple algebraic group. Given $g\in G$ and an invariant function $f\in \C[G]^G$, the corresponding element $F_g$ lies in the Lie algebra $\Lie(Z_G(g))$ of the centraliser $Z_G(g)$.
\end{proposition}
\begin{proof}
Recall that 
\[\Lie(Z_G(g)) = \{\alpha\in \frakg\,|\, \Ad_g(\alpha) = \alpha\} = \ker(\Ad_g-I:\frakg\to \frakg).\]
It suffices to show that 
\[F_g\in \bigl(\im(\Ad_g-I)\bigr)^\perp= \ker(\Ad_g - I),\]
where $(-)^\perp$ denotes the orthogonal complement with respect to the form $\kappa$.
To this end, we compute
\[\kappa(F_g, \Ad_g\alpha - \alpha) = \kappa(F_g, \Ad_g\alpha) -\kappa(F_g, \alpha).\]
We are done if we show the right hand side is zero.
We have
\begin{align*}
\kappa(F_g,\Ad_g\alpha) = \frac{d}{dt}\bigg|_{t=0}f(g\exp(t\Ad_g(\alpha)))&=\frac{d}{dt}\bigg|_{t=0}f(gg\exp(t\alpha)g^{-1})\\
&=\frac{d}{dt}\bigg|_{t=0}f(g\exp(t\alpha)) = \kappa(F_g,\alpha),
\end{align*}
and the statement follows.
\end{proof}

\begin{example}
Let $G=\SL_n(\C)$, $f=\tr$, and $g\in G$. We will show that $F_{g} = g-\frac{\tr(g)}{n}I\in \Lie(Z_G(g))$. Instead of using the Killing form, in this example we use the non-degenerate invariant bilinear form $\kappa(X,Y)=\tr(XY)$. We have $\kappa(F_g,X ) = \tr(F_g X)$ and, on the other hand, we have
\[\kappa(F_g, X) = \frac{d}{dt}\bigg|_{t=0}f(g\exp(tX)) = \tr(gX),\quad \text{for all}\,\,X\in \sl_n(\C).\]
Thus, 
\[0=\tr(F_g X)-\tr(gX) = \tr((F_g - g)X).\]
Hence we conclude that $F_g - g$ is orthogonal to $\sl_n(\C)$. Using the standard decomposition $\gl_n(\C)  = \sl_n(\C) \oplus \C I_n$, the orthogonal complement of $\sl_n(\C)$ with respect to the trace form is $\C I_n$. Thus,
\[F_g-g=\lambda I_n\]
for some $\lambda\in \C$. Taking traces on both sides, one finds that $\lambda = -\frac{\tr(g)}{n}$. Thus, 
\[F_g = g-\frac{\tr(g)}{n}I_n.\]
\end{example}

Since $G$ is semisimple, its algebra of invariant functions is a polynomial algebra
\[\C[G]^G = \C[\chi_1,\dots, \chi_r]\]
where $r = \rank(G)$ and the $\chi_i$ are the characters of the fundamental representations (cf. \cite{Ste65}). For instance, when $G=\SL_n(\C)$, the generators are the coefficients of the characteristic polynomial. We end this section with a proposition due to Steinberg which will be useful for defining a basis for the Lie algebras of the centralisers of regular elements:
\begin{theorem}{\cite[Theorem 8.1]{Ste65}}\label{thm:Steinberg}
Let $G$ be a semisimple simply connected algebraic group. An element $x\in G$ is regular if and only if the differentials of the fundamental characters $d\chi_1,\dots,d\chi_r$ are linearly independent at $x$.
\end{theorem}

\section{Regular commuting variety}\label{c:regularcommutingvariety}
In this section, we prove that the regular commuting variety lies inside the smooth locus of the commuting variety and that it has LSAP, and deduce the Oka-$1$ property of the smooth locus.
\subsection{Smoothness of the regular commuting variety}\label{c:Smoothness}

Let $G$ be a simply connected semisimple algebraic group and let $R=R(G)$ denote the commuting variety defined as
\[R = \{(x,y)\in G\times G\,|\, xy=yx\}.\]
There is a projection map from $R$ to $G$: 
\[\pi:R\longrightarrow G,\quad (x,y)\longmapsto x.\]
 Consider the inverse image under $\pi$ of the dense open subset of regular elements $G_{\rm reg}\subset G$
 \[R_{\mathrm{reg}} := \pi^{-1}(G_{\rm reg}) =  \{(x,y)\in R\,|\, x\in G_{\mathrm{reg}}\}.\]
We call the open subset $R_{\mathrm{reg}}\subset R$ the \emph{regular commuting variety}. Furthermore, $R_{\rm reg}$ is irreducible since it is an open subset of the irreducible variety $R$ \cite{Ric79}. Now we have a surjective map $\pi|_{R_{\rm reg}}:R_{\rm reg}\to G_{\rm reg}$ between irreducible varieties whose fibre over each point $x\in G_{\rm reg}$ has dimension equal to $\rank(G)$. By \cite[\S 6.3, Theorem 1.25]{Sha13} it follows that 
\[\dim R = \dim R_{\rm reg}=\dim(G)+\rank(G).\]
We identify $T_xG$ and $T_yG$ with $\frakg=\Lie(G)$ by right translation, and hence $T_{(x,y)}(G\times G)$ with $\frakg^2$.
Consider the set
\[K_{(x,y)}= \{(\alpha,\beta)\in \frakg^2\,|\, \Ad_x(\beta)-\beta=\Ad_y(\alpha)-\alpha\}.\]
Then differentiating the commuting equation $xy=yx$ gives $T_{(x,y)}R\subseteq K_{(x,y)}$.
We have the following proposition:
 \begin{proposition}\label{prop:submersion}
For every $(x,y)\in R_{\rm reg}$ the projection  
 \[p:K_{(x,y)} \longrightarrow T_xG\]
 onto the first coordinate is surjective. Moreover, the point $(x,y)$ is smooth.
 \end{proposition}
 
 \begin{proof}
To show that $p$ is surjective, we need, for a given $\alpha\in \frakg$, to show that there is $\beta\in \frakg$ with
\[\Ad_x(\beta) - \beta = \Ad_y(\alpha) - \alpha.\]
Consider the linear map
\[T_x:\frakg \longrightarrow \frakg, \quad \beta \longmapsto \Ad_x(\beta) - \beta.\]  
It suffices to show that $\Ad_y(\alpha) - \alpha$ lies in $\im(T_x)$.

Let $\kappa$ denote the Killing form on $\frakg$, which is non-degenerate because $\frakg$ is semisimple. By basic linear algebra,
\[(\im T_x)^\perp = \ker(T_x^*),\]  
where $(\im T_x)^\perp$ denotes the orthogonal complement of $\im(T_x)$ and $T_x^*$ the adjoint with respect to $\kappa$.

We claim that $T_x^* = T_{x^{-1}}$. Indeed, using the $\Ad$-invariance of $\kappa$, for $\beta,\gamma\in \frakg$ we compute
\[\kappa(T_x(\beta), \gamma) = \kappa(\Ad_x(\beta) - \beta, \gamma) = \kappa(\beta, \Ad_{x^{-1}}(\gamma) - \gamma) = \kappa(\beta, T_{x^{-1}}(\gamma)),\]
so $T_x^* = T_{x^{-1}}$. Hence,
\[\im(T_x) = (\ker T_{x^{-1}})^\perp = (\ker T_x)^\perp,\]
where the last equality holds because $\Ad_x(\beta) = \beta$ if and only if $\Ad_{x^{-1}}(\beta) = \beta$.

Therefore, it suffices to show that
\[\kappa(z, \Ad_y(\alpha) - \alpha) = 0 \quad \text{for all } z\in \ker(T_x).\]  
But for $z\in \ker(T_x)$, using $\Ad$-invariance of $\kappa$,
\[\kappa(z, \Ad_y(\alpha) - \alpha) = \kappa(z, \Ad_y(\alpha)) - \kappa(z, \alpha) = \kappa(\Ad_{y^{-1}}(z) - z, \alpha).\]
Finally, since $x\in G_{\mathrm{reg}}$, its centraliser $Z_G(x)$ in $G$ is abelian \cite[Theorem~D]{Kur78}. As $y$ commutes with $x$, we have $y\in Z_G(x)$, so $y$ acts trivially on $\Lie(Z_G(x)) = \ker(T_x)$. Hence, $\Ad_{y^{-1}}(z) = z$ for all $z\in \ker(T_x)$, giving
\[\kappa(\Ad_{y^{-1}}(z) - z, \alpha) = 0.\]  
This proves that $\Ad_y(\alpha) - \alpha \in \im(T_x)$, so $p$ is surjective. Consequently, 
\[\dim (K_{(x,y)}) = \dim G + \dim(\ker(p)) = \dim G +\rank(G).\]
Moreover, for $(x,y)\in R_{\rm reg}$ we have
\[\dim R\leq \dim(T_{(x,y)}R)\leq \dim(K_{(x,y)})=\dim R,\]
 so $(x,y)$ is a smooth point in $R$ and thus $R_{\rm reg}$ is contained in the smooth locus of $R$.
\end{proof}

\subsection{Local parameterisation of the regular commuting variety}\label{c:localparams}
In this section, we give local parameterisations for the open subset $R_{\rm reg}$ of the commuting variety of a simply connected almost simple algebraic group $G$. The more general simply connected semisimple case will follow by taking products over the almost simple factors. For $x\in G$, we write $\frakz(x):=\Lie(Z_G(x))$ for the Lie algebra of its centraliser.

Let $F_1,\dots,F_r:G_{\rm reg}\to \frakg$ be holomorphic maps associated to a set of generators for the invariant functions on $G$, so that $F_{x,i}:=F_i(x)$, for $i=1,\dots,r$, form a basis for $\frakz(x)$ (cf. \S~\ref{s:differentialsofinvariantfunctions}).

\begin{lemma}\label{lem:bihol}
Let $G$ be a simply connected almost simple algebraic group of rank $r$. For a fixed integer $n$ and $\varepsilon\in Z(G)$, the map
\[
\Phi_{(n,\varepsilon)}:G_{\mathrm{reg}}\times \C^r \longrightarrow R_{\mathrm{reg}},
\qquad
(x,\lambda_1,\dots,\lambda_r)\longmapsto \bigl(x,\varepsilon \,x^n\exp\bigl(\sum_{i=1}^r\lambda_i F_{x,i}\bigr)\bigr),
\]
is a local biholomorphism onto its image.
\end{lemma}

\begin{proof}
Since
\[ \dim R_{\rm reg}=\dim G+\rank(G) = \dim(G_{\rm reg}\times \mathbb C^r),\]
by the holomorphic inverse function theorem, it suffices to show that
the differential of \(\Phi_{(n,\varepsilon)}\) is injective at every point of $G_{\rm reg}\times \C^r$.

Fix \(\varepsilon\in Z(G)\) and \(n\in\mathbb Z\). For ease of notation, write
\[ \Phi=\Phi_{(n,\varepsilon)}\quad \text{and}\quad (x,\lambda):=(x,\lambda_1,\dots,\lambda_r).\]
Let $(\alpha,\beta_1,\dots,\beta_r)\in T_xG_{\rm reg}\oplus \C^r$ and suppose that the differential at $(x,\lambda)\in G_{\rm reg}\times \C^r$ satisfies
\[
        d\Phi_{(x,\lambda)}(\alpha,\beta_1,\dots,\beta_r)=0.
\]
The restriction of $\Phi$ to the first coordinate is the identity map on \(G_{\rm reg}\), and hence so is its differential. Then, by assumption, we must have $\alpha=0$. In the second coordinate, all terms involving $d(x^n)_x(\alpha)$ and $(dF_{i})_x(\alpha)$ vanish. Thus the second component of $d\Phi_{(x,\lambda)}(\alpha,\beta_1,\dots,\beta_r)$ reduces to
\[\varepsilon x^n d\exp_{A(x,\lambda)}\left(\sum_{i=1}^r \beta_iF_{x,i} \right) = 0.\]
We now compute the differential $d\exp_{A(x,\lambda)}\left(\sum_{i=1}^r \beta_iF_{x,i} \right)$. Write
\[A(x,\lambda):=\sum_{i=1}^r\lambda_iF_{x,i}
\in \frakz(x)\quad \text{and}\quad   B:=\sum_{i=1}^r\beta_iF_{x,i}\]
and observe that $A(x,\lambda)$ and $B$ both lie in $\frakz(x)$ (cf. Proposition~\ref{prop:differentialinvariantfunction}). Since $x\in G_{\rm reg}$ and $G$ is simply connected, the Lie algebra $\frakz(x)$ is abelian (\cite[Remark 3.13]{Ste65} or \cite[Theorem D]{Kur78}). In particular, $A(x,\lambda)$ and $B$ commute, so we have
\[\exp(A(x,\lambda)+tB)=\exp(A(x,\lambda))\exp(tB).\]
So
\[d\exp_{A(x,\lambda)}(B)=\left.\frac{d}{dt}\right|_{t=0}\exp(A(x,\lambda)+tB)= \exp(A(x,\lambda))B.\]
In particular, we get
\[
\varepsilon x^n\exp(A(x,\lambda))B=0.
\]
Since \(\varepsilon x^n\exp(A(x,\lambda))\) is invertible,
it follows that $ B=0$, that is,
\[\sum_{i=1}^r\beta_iF_{x,i}=0.\]
By Theorem~\ref{thm:Steinberg}, $F_{x,1},\dots,F_{x,r}$ are linearly independent in \(\mathfrak z(x)\), so $\beta_i=0$ for all $i=1,\dots,r$.
Thus
\[(\alpha,\beta_1,\dots,\beta_r)=0,\]
and hence $d\Phi_{(x,\lambda)}$ is injective.
\end{proof} 

\begin{lemma}\label{lem:Phi-cover-Rreg}
Let $G$ be a simply connected almost simple algebraic group. Then
\[R_{\rm reg} = \bigcup_{(n,\varepsilon)\in \Z\times Z(G)}\im(\Phi_{(n,\varepsilon)}).\]
\end{lemma}

\begin{proof}
Let $(x,y)\in R_{\rm reg}$. Then $x\in G_{\rm reg}$ and $y\in Z_G(x)$. By Lemma~\ref{lem:regular-centraliser-components}, there exist
$n\in\mathbb Z$, $\varepsilon\in Z(G)$, and $h\in Z_G(x)^\circ $ such that $y=x^n\varepsilon h$. We will show that $h$ lies in the image of $\exp:\frakz(x)\to Z_G(x)^\circ$. Since $Z_G(x)^\circ$ is a connected commutative linear algebraic group, it is isomorphic to a product $(\C^\times\!)^s\times U,$ where $U$ is a connected commutative unipotent group \cite[Section 1.6]{Car93}. The exponential maps from the Lie algebras $\Lie(\C^\times\!)$ and $\Lie(U)$ to the respective groups $\C^\times$ and $U$ are surjective. Hence, the exponential map
\[\exp:\mathfrak z(x)\longrightarrow Z_G(x)^\circ\]
is surjective. Therefore there exists an $ \alpha\in\frakz(x)$ with $h=\exp(\alpha)$. Since $F_{x,1},\dots,F_{x,r}$ form a basis of $\frakz(x)$, we can write
\[\alpha=\sum_{i=1}^r \lambda_iF_{x,i}\]
for some $\lambda_1,\dots,\lambda_r\in\C$ and, in particular,
\[y=x^n\varepsilon\exp\left(\sum_{i=1}^r \lambda_iF_{x,i}\right).\]
Hence
\[(x,y) = \Phi_{(n,\varepsilon)}(x,\lambda_1,\dots,\lambda_r).\]
Thus, the images of the maps $\Phi_{(n,\varepsilon)}$ cover $R_{\rm reg}$.
\end{proof}
We will require a special case of the following lemma which says that a complex manifold covered by images of local biholomorphisms from LSAP manifolds has LSAP.

\begin{lemma}\label{lem:LSAP-local-bihol}
Let $(X_\alpha)_{\alpha\in I}$ be a family of complex manifolds, each having LSAP, and let $Y$ be another complex manifold. Suppose further for each $\alpha$ we have a collection $\{\Phi_{\alpha,\beta}: X_{\alpha}\to Y\}_{\beta\in J}$ of local biholomorphisms indexed by a set $J$ such that 
\[\bigcup_{(\alpha,\beta)\in I\times J}\im(\Phi_{\alpha,\beta}) =Y.\]
Then $Y$ has the LSAP property.
\end{lemma}

\begin{proof}
Let $y_0\in Y$. By the assumption in the lemma, there are $\alpha\in I$, $\beta\in J$, $x_0\in X_\alpha$ and neighbourhoods $U\subset X_\alpha$ of $x_0$ and $W\subset Y$ of $y_0$ such that 
\[\Phi_{\alpha,\beta}:U\longrightarrow W\] 
is a biholomorphism. Since \(X_\alpha\) has LSAP, after possibly shrinking \(U\), we may assume that \(U\) is an LSAP neighbourhood of \(x_0\). We claim that $W$ is an LSAP neighbourhood of $y_0$. Let
\(D\subset D'\subset \C\) be a special pair of compact discs, and let
\[
H:D\times \BB^N\longrightarrow W
\]
be a dominating holomorphic spray. Composing with \(\Phi_{\alpha,\beta}^{-1}\), we
obtain a dominating holomorphic spray
\[
\widetilde H:=\Phi_{\alpha,\beta}^{-1}\circ H:
D\times \BB^N\longrightarrow U\subset X_\alpha.
\]
By the LSAP property for $X_\alpha$, there exists $r=r(\widetilde H)\in(0,1)$
such that \(\widetilde H\) can be approximated uniformly on $D\times r\overline{\BB}^N$ by holomorphic maps
\[
\widetilde G:D'\times r\overline{\BB}^N\longrightarrow X_\alpha.
\]
Composing with \(\Phi_{\alpha,\beta}\), we obtain holomorphic maps
\[
G:=\Phi_{\alpha,\beta}\circ \widetilde G:
D'\times r\overline{\BB}^N\longrightarrow Y.
\]
Since $\Phi_{\alpha,\beta}$ is holomorphic, these maps approximate $H$ uniformly on $D\times r\overline{\BB}^N$. Therefore $W$ is an LSAP neighbourhood of $y_0$, so $Y$ has the LSAP property at $y_0$. Since $y_0\in Y$ was arbitrary, $Y$ has the LSAP property.
\end{proof}

\subsection{Oka-$1$ theorem}\label{c:Oka1theorem}

\begin{theorem}\label{thm:RregLSAP}
Let $G$ be a connected simply connected semisimple algebraic group. Then $R_{\rm reg}$ has the LSAP property.
\end{theorem}

\begin{proof}
Write $G=G_1\cdots G_m$ as the product of its almost simple factors. Since $G$ is simply connected, the multiplication map
\[
G_1\times\cdots\times G_m\longrightarrow G
\]
is an isomorphism \cite[\S 24a]{Mil17}. Under this isomorphism,
\[
G_{\rm reg}\cong\prod_{j=1}^m (G_j)_{\rm reg}
\qquad\text{and}\qquad
R_{\rm reg}(G)\cong\prod_{j=1}^m R_{\rm reg}(G_j).
\]
Applying Lemmas~\ref{lem:bihol} and~\ref{lem:Phi-cover-Rreg} to each almost simple factor and taking products, we obtain a family of local biholomorphisms from $G_{\rm reg}\times\C^r$ whose images cover $R_{\rm reg}(G)$.

By Lemma~\ref{lem:LSAP-local-bihol}, it remains to show that $G_{\rm reg}\times \C^r$ has LSAP. Since finite products of LSAP manifolds have LSAP and $\C^r$ has LSAP, it is enough to show that $G_{\rm reg}$ does. Indeed, in \cite[Proposition 1.1]{AZK12} it is shown that semisimple algebraic groups are flexible and in \cite[Theorem 1.1]{FKZ16} that complements of flexible affine varieties of a closed subvariety of codimension $\geq 2$ remain flexible. Finally, Steinberg \cite[Theorem 1.3 (a) and (b)]{Ste65} showed that $G_{\rm reg}$ is the complement of a closed subvariety of $G$ of codimension $\geq 3$. Thus, $G_{\rm reg}$ is flexible, hence has the LSAP property.
 \end{proof}
\begin{remark}
Recall that we defined $R_{\rm reg}$ as the preimage of $G_{\rm reg}$ by the projection map $\pi:G\times G\to G$ on the first component. However, the choice of projection map to the first component was arbitrary. If we denote the two projection maps $\pi_i$ for $i=1,2$, then $R_{\rm reg} = \pi^{-1}_{1}|_{R(G)}(G_{\rm reg})$. Define
\[R_{\rm reg}^{(i)} := \pi_i^{-1}|_{R(G)}(G_{\rm reg}).\]
By symmetry the conclusion of Theorem~\ref{thm:RregLSAP} also holds for 
\[R_{\rm reg}^{(2)} = \{(x,y)\in R\,|\, y\in G_{\rm reg}\}.\]
Since the LSAP property is a Hausdorff local property \cite[Proposition 5.5]{FL25}, we have that the union $R^{(1)}_{\rm reg}\cup R^{(2)}_{\rm reg}$ has LSAP. 
\end{remark}
\begin{corollary}\label{cor:mainresult} Let $G$ be a simply connected semisimple algebraic group. The smooth locus $R(G)_{\rm sm}$ of the commuting variety of $G$ is an Oka-$1$ manifold.
\end{corollary}
\begin{proof}
Set $n= \dim_\C(R_{\rm sm})$. The complement $E =R_{\rm sm}\setminus R_{\rm reg}$ is a proper analytic subset of $R_{\rm sm}$ and thus has complex dimension $\dim_\C(E)\leq n-1$. Thus, the Hausdorff measure $\calH^{2n-1}_{R_{\rm sm},g}(E)$ of $E$ is zero. By Theorem~\ref{thm:RregLSAP}, $R_{\rm reg}$ has LSAP and hence $R_{\rm sm}$ is Oka-$1$ by \cite[Proposition 7.13  (b)]{AF25}.
\end{proof}

\begin{remark}
\begin{enumerate}
\item There are many instances where the assumption of simple connectedness of $G$ was used. One of the main reasons we require simple connectedness is Steinberg's result that, for $G$ simply connected, $F_{x,1},\dots,F_{x,r}\in \frakz(x)$ are linearly independent if and only if $x\in G_{\rm reg}$. Other applications of simple connectedness are that the centralisers of semisimple elements are connected and that the centralisers of regular elements are abelian. 
\item Within the classification of reductive algebraic groups, the simply connected almost simple algebraic groups are $\SL_{n+1}(\C)$, $\Spin_{2n+1}(\C)$, $\Sp_{2n}(\C)$, $\Spin_{2n}(\C)$ together with the simply connected forms of the exceptional types $G_2,$ $F_4$, $E_6$, $E_7$, $E_8$ \cite[Table~9.2, p.~72]{MT11}. Every simply connected semisimple algebraic group is a finite direct product of groups from this list \cite[\S 24a]{Mil17}.
\end{enumerate}
\end{remark}

\section{Strongly commuting varieties}\label{c:stronglycommutingvarieties}
The motivation for defining strongly commuting varieties is the following. Suppose that we are interested in commuting elements in $\PSL_n(\C)$. Two elements $[x],[y]\in \PSL_n(\C)$ commute if and only if their representatives $x,y\in \SL_n(\C)$ commute up to an element of the centre $\omega\in Z(\SL_n(\C))$, i.e., $xy=\omega yx.$
Hence, we may impose a stronger condition on the elements $[x],[y]\in \PSL_n(\C)$ to commute by requiring that the representatives $x,y\in \SL_n(\C)$ commute in $\SL_n(\C)$. In this case, we will say that $[x]$ and $[y]$ strongly commute. We give four equivalent descriptions of which two are topological and two are algebraic in nature. The algebraic formulations are inspired by \cite[Section 3]{Ric79}.

We begin with some notation. Let $G$ be a connected semisimple algebraic group and let $q:\widetilde{G}\to G$ be its universal covering, and let $\Gamma:=\ker q$. Since $\widetilde{G}$ is connected and $\Gamma$ is a finite normal subgroup, it follows that $\Gamma\subset Z(\widetilde{G})$. Let $(x,y)\in R(G)$ and $(X,Y)\in \widetilde{G}^2$ be any lift. Define the set-theoretic map
\[\kappa: R(G)\longrightarrow \Gamma,\quad (x,y)\longmapsto [X,Y].\]
\begin{proposition}\label{prop:kappa}
The map $\kappa$ is well-defined and is a morphism of varieties.
\end{proposition}
Before proving this proposition, for each $\gamma\in \Gamma$, we introduce
\[R_\gamma(\widetilde{G}) = \{(X,Y)\in \widetilde{G}^2\,|\,[X,Y] = \gamma\}.\]
Then, we have
\[(q\times q)^{-1}(R(G)) =\bigsqcup_{\gamma\in \Gamma}R_\gamma(\widetilde{G}).\]
Indeed, for any $(X,Y)\in R_\gamma(\widetilde{G})$ the images $q(X)$ and $q(Y)$ commute. Conversely, for any $(X,Y)\in (q\times q)^{-1}(x,y)$ with $(x,y)\in R(G)$, we have that $[X,Y]=\gamma$. Applying $q$ to both sides shows that $\gamma\in \Gamma$. Now define the following subset of $R(G)$
\[R_\gamma(G) := \kappa^{-1}(\gamma).\]
We think of $R_\gamma(G)$ as the subset of $R(G)$ whose lifts to $\widetilde{G}$ commute up to the element $\gamma\in \Gamma$. It follows immediately that $R_\gamma(G) = (q\times q)(R_\gamma(\widetilde{G}))$ and $R_\gamma(G)\cap R_\delta(G) = \emptyset$ if $\gamma\neq \delta$. 
\begin{proof}[Proof of Proposition~\ref{prop:kappa}]
Since $q\times q:\widetilde{G}\times \widetilde{G}\to G\times G$ is finite and hence closed, the $R_\gamma(G)=(q\times q)(R_\gamma(\widetilde{G}))$ are closed in $R(G)$, as the $R_\gamma(\widetilde{G})$ are closed. Since there are finitely many $\gamma\in \Gamma$, the $R_\gamma(G)$ are pairwise disjoint and their union covers $R(G)$. It follows that the $R_\gamma(G)$ are also open. The result of the proposition follows from the following claim: A set-theoretic map $f:X\to \Gamma$ between an affine variety $X\subset\AA^n$ and a finite variety $\Gamma\subset \AA^m$ is induced by a morphism of varieties if and only if the fibres are open and closed. The forward direction is obvious. The converse is an application of \cite[II, Ex 2.19]{Har77}. 
\end{proof}
\begin{definition}\label{def:str}
Two commuting elements $a,b\in G$ \emph{strongly commute} if 
\[(a,b)\in (q\times q)(R(\widetilde{G})).\]
We define the \emph{strongly commuting variety} of $G$ as 
\[R_{\rm str}(G) := (q\times q) (R(\widetilde{G}))=R_1(G).\]
\end{definition}
By Richardson's theorem, $R(\widetilde{G})$ is irreducible because $\widetilde{G}$ is semisimple and simply connected. Therefore, the strongly commuting variety is the image of the irreducible variety, hence is irreducible and connected and contains $(1,1)$. Thus, it is the identity component of $R(G)$. The topology of the identity component of $R(G)$ has been studied in \cite{Bai07,GPS12,PS13}.

We now state another topological description.

\begin{definition}\label{def:top}
Two commuting elements $a,b\in G$ \emph{topologically strongly commute} if for some paths $\alpha,\beta :[0,1]\to G,$ with 
\begin{equation}\label{eq:path}
\alpha(0)=\beta(0)=1\,\, \text{and}\,\,
 \alpha(1)=a\,\, \text{and}\,\, \beta(1)=b
 \end{equation}
the product paths $\alpha(t)\beta(t)$ and $\beta(t)\alpha(t)$ in $G$ are homotopic relative to the endpoints. Denote by $R_{\rm top}(G)$ the collection of pairs of elements which topologically strongly commute.
\end{definition}

The next definition is algebraic and comes from \cite[\S 3]{Ric79}.
\begin{definition}\label{def:alg}
Let $\calE'(G)$ denote the pairs of elements $(x,y)\in G^2$ such that there exists a maximal torus $T$ containing both $x$ and $y$. Let $\calE(G)$ denote the Zariski closure of $\calE'(G)$.
\end{definition}

We now introduce a fourth definition which is also inspired by \cite[\S 3]{Ric79}:
\begin{definition}
Consider the subset 
\[S'(G) = \{(x,y)\in R(G)\,|\, \exists\, z\in Z(G)\,\, \text{such that}\,\, zy\in Z_G(x)^\circ\}\]
of $R(G)$. Denote the Zariski closure of $S'(G)$ by $S(G)$.
\end{definition}

\begin{theorem}
Let $G$ be a semisimple algebraic group. We have
\[R_{\rm str}(G)=R_{\rm top}(G)=\calE(G)=S(G).\]
\end{theorem}

We prove this in the following order: $R_{\rm str}(G)=R_{\rm top}(G)$ then $R_{\rm str}(G)=\calE(G)$, and lastly $\calE(G)=S(G)$. 

\medskip

\begin{proposition}
The equality $R_{\rm str}(G) = R_{\rm top}(G)$ holds.
\end{proposition}

\begin{proof}
Let $(a,b)\in R_{\rm str}(G)$, i.e., $a,b\in G$ strongly commute.  We will show that $(a,b)\in R_{\rm top}(G)$, i.e., $a$ and $b$ topologically strongly commute. Let $\alpha,\beta: [0,1]\to G$ be paths with 
\[
\alpha(0)=\beta(0)=1\,\, \text{and}\,\,
 \alpha(1)=a\,\, \text{and}\,\, \beta(1)=b.
\]
We must show that the product paths $\alpha\beta$ and $\beta\alpha$ are homotopic. By assumption, there are lifts $\tilde{a},\tilde{b}\in \tilde{G}$ such that
\[\tilde{a}\tilde{b} = \tilde{b}\tilde{a}.\]
Let $\tilde{\alpha},\tilde{\beta}:[0,1]\to \tilde{G}$ be the unique lifts of $\alpha$ and $\beta$ such that 
\[\tilde{\alpha}(0)=1_{\tilde{G}}=\tilde{\beta}(0).\]
Then $\tilde{\alpha}(1)$ and $\tilde{a}$ differ by an element $\hat{u} $ in the kernel $\ker(q) $. Specifically, if $\tilde{\alpha}(1)=\tilde{u}$ then $q(\tilde{a}\tilde{u}^{-1})= 1_G$. Similarly, for $\tilde{\beta}(1)=\tilde{v}$, we have $\tilde{b}\tilde{v}^{-1}\in\ker q$. Define $\hat{u}:=\tilde{a}\tilde{u}^{-1}$ and $\hat{v}:=\tilde{b}\tilde{v}^{-1}$. Since finite normal subgroups of connected algebraic groups are central, it follows that $\hat{u},\hat{v}\in \ker(q)$ are central elements in $\tilde{G}$. Moreover, from the assumption that $\tilde{a}\tilde{b}=\tilde{b}\tilde{a}$ it follows that
\[\tilde{\alpha}(1)\tilde{\beta}(1) = \tilde{\beta}(1)\tilde{\alpha}(1).\]
Thus, the endpoints of the lifts $\tilde{\alpha}$ and $\tilde{\beta}$ commute so the products of the lifts $\tilde{\alpha}\tilde{\beta}$ and $\tilde{\beta}\tilde{\alpha}$
have the same endpoints and therefore are homotopic by some homotopy $\tilde{H}$ as $\tilde{G}$ is simply connected. Then $q\circ \tilde{H}$ gives a homotopy between $\alpha\beta$ and $\beta\alpha$. Thus, $a$ and $b$ topologically strongly commute.

Conversely, suppose $a$ and $b$ topologically strongly commute. We are required to show that there exist $\tilde{a},\tilde{b}\in \tilde{G}$ such that $\tilde{a}\tilde{b} = \tilde{b}\tilde{a}$ and $q(\tilde{a})=a$ and $q(\tilde{b})=b$. 

The paths $\alpha,\beta:[0,1]\to G$ have unique lifts $\tilde{\alpha},\tilde{\beta}:[0,1]\to \tilde{G}$ with start point $1_{\tilde{G}}$. The product paths $\tilde{\alpha}\tilde{\beta}$ and $\tilde{\beta}\tilde{\alpha}$ satisfy
\[q(\tilde{\alpha}\tilde{\beta}) = \alpha\beta\quad \text{and}\quad q(\tilde{\beta}\tilde{\alpha})=\beta\alpha\]
because $q$ is a homomorphism of groups so $\tilde{\alpha}\tilde{\beta}$ and $\tilde{\beta}\tilde{\alpha}$ are the unique lifts for the paths $\alpha\beta$ and $\beta\alpha$. Since $\tilde{\alpha}(0)\tilde{\beta}(0)=1_{\tilde{G}} = \tilde{\beta}(0)\tilde{\alpha}(0)$ by the covering homotopy lemma, 
\[\tilde{\alpha}(t)\tilde{\beta}(t) \quad \text{and} \quad \tilde{\beta}(t)\tilde{\alpha}(t)\]
have the same endpoints and are homotopic.  In particular, 
\[\tilde{\alpha}(1)\tilde{\beta}(1)=\tilde{\beta}(1)\tilde{\alpha}(1).\]
Thus, $(\tilde{\alpha}(1),\tilde{\beta}(1))\in R(\tilde{G})$ and 
\[(q\times q)(\tilde{\alpha}(1),\tilde{\beta}(1))=(a,b)\in R(G).\qedhere\] 
\end{proof}

\begin{proposition}
The equality $R_{\rm str}(G) = \calE(G)$ holds.
\end{proposition}

\begin{proof}
If two elements $\tilde{a},\tilde{b}$ are contained in a maximal torus $\tilde{T}\subset \tilde{G}$ then their images $q(\tilde{a})$ and $q(\tilde{b})$ are contained in the maximal torus $T:=q(\widetilde{T})$ of $G$ \cite[Proposition 11.4]{Bor91}. Conversely, if $T\subset G$ is a maximal torus, then $\widetilde{T}=q^{-1}(T)^\circ $ is a maximal torus of $\widetilde{G}$ with $q(\widetilde{T})=T$. Hence any pair in $T^2$ has lifts in $\widetilde{T}^2$. Since the map $q\times q$ is a finite map, it is closed and, in particular, the map commutes with taking Zariski closures. With these facts at hand we have
\[\calE(G) = \overline{\calE'(G)}=\overline{(q\times q)(\calE'(\tilde{G}))}=(q\times q)(\calE(\tilde{G})) = (q\times q)(R(\tilde{G}))=R_{\rm str}(G).\]
Note that the fourth equality is \cite[Theorem C]{Ric79}.
\end{proof}

\begin{proposition}
The equality $\calE(G) = S(G)$ holds.
\end{proposition}

\begin{proof}
The inclusion $S(G)\subset \calE(G)$ follows from \cite[\S 3, Theorem B]{Ric79}. For the other inclusion, suppose $x$ and $y$ are contained in the maximal torus $T$ so that $(x,y)\in \calE'(G)$. Then $y\in Z_G(x)$ and $T\subset Z_G(x)^\circ$, thus, $y\in Z_G(x)^\circ$. In particular, $(x,y)\in S'(G)$. Thus, $\calE(G)\subset S(G)$.
\end{proof}

\subsection{Oka-$1$ for strongly commuting varieties}\label{c:Oka1theoremstronglycommutingvarieties}
In this section we prove that the smooth locus of the strongly commuting variety of a semisimple algebraic group is Oka-$1$. Let $G$ be a semisimple algebraic group and $R_{\rm str}(G)$ be its strongly commuting variety. The open subset of $R_{\rm str}(G)$ defined by
\[R_{\rm str,reg}(G):=\{(x,y)\in R_{\rm str}(G)\,|\, x\in G_{\rm reg}\}\]
is called the \emph{regular strongly commuting variety}.
\begin{theorem}\label{thm:strongOka1}
Let $G$ be a connected semisimple algebraic group. Then the regular strongly commuting variety has LSAP. Consequently, the smooth locus $R_{\rm str}(G)_{\rm sm}$ is an Oka-$1$ manifold.
\end{theorem}

\begin{proof}
Since $q:\widetilde{G}\to G$ is a central isogeny, it is a finite holomorphic covering map hence a local biholomorphism. Moreover, regular elements in $\widetilde{G}$ descend to regular elements in $G$ via $q$ \cite[Proposition 2.3 (2)]{Kur78}. It follows that the map
\[q\times q:R_{\rm reg}(\tilde{G})\longrightarrow R_{\rm str, reg}(G)\]
is a surjective local biholomorphism. Since $R_{\rm reg}(\widetilde{G})$ is smooth and has LSAP by Theorem~\ref{thm:RregLSAP}, it follows that $R_{\rm str, reg}(G)$ is smooth and also has LSAP. Set $n=\dim_\C(R_{\rm str}(G)_{\rm sm})$. The complement $E =R_{\rm str}(G)_{\rm sm}\setminus R_{\rm str, reg}(G)$ is a proper analytic subset of complex dimension $\dim_\C(E)\leq n-1$. Thus the Hausdorff measure $\calH^{2n-1}_{R_{\rm str}(G)_{\rm sm},g}(E)$ of $E$ is zero. Thus, $R_{\rm str}(G)_{\rm sm}$ is Oka-$1$ by \cite[Proposition 7.13  (b)]{AF25}.
\end{proof}

\subsection{Oka-$1$ for $\PSL_p(\C)$}\label{c:Oka1PSLp}
In Theorem~\ref{thm:strongOka1} we showed that the smooth locus of the strongly commuting variety of a semisimple algebraic group is Oka-$1$. In this section, we prove an Oka-$1$ result for the full commuting variety of the semisimple (but not simply connected) group $\PSL_p(\C)$ when $p$ is prime. This will involve the study of \emph{almost commuting elements} which have previously been studied in \cite{BFM02} and \cite{KS00} in the case of compact Lie groups. Here, we adapt some of their results to the setting of semisimple algebraic groups.

\medskip 

We return to the case where $G$ is a semisimple algebraic group and let $q:\widetilde{G}\to G$ be its universal covering and let $\Gamma=\ker(q)$. Then
\[(q\times q)^{-1}(R(G)) = \bigsqcup_{\gamma\in \Gamma}R_\gamma(\widetilde{G}),\]
where 
\[R_{\gamma}(\widetilde{G}) = \{(X,Y)\in \widetilde{G}^2\,|\, [X,Y]=\gamma\}\]
and we have identified the strongly commuting variety with the image $(q\times q)(R_1(\widetilde{G}))$. It is an interesting problem to describe the components $R_\gamma(\widetilde{G})$, i.e., collections of elements commuting up to a central element $\gamma\in \Gamma$. Following Serre \cite{Ser05}, we define the notion of irreducible pairs:  
\begin{definition}
A subgroup $H$ of $G$ is said to be \emph{irreducible} if it is not contained in any proper parabolic subgroup of $G$. We say that a pair $(x,y)\in G\times G$ is \emph{irreducible} if $\overline{\langle x,y\rangle}$ is an irreducible subgroup of $G$. Here $\overline{\langle x,y\rangle}\subset G$ denotes the Zariski closure of the group generated by elements $x,y\in G$.
\end{definition}

Note that irreducibility of $(x,y)$ implies that $H:= \overline{\langle x,y\rangle}$ is reductive. Indeed, the unipotent radical $R_u(H)$ of $H$ is normalised by $H$, and hence $H\subset N_G(R_u(H))$. By \cite[Cor 30.3]{Hum75}, $N_G(R_u(H))$ is contained in a parabolic subgroup of $G$, so, by irreducibility of $H$, it follows that $R_u(H)=\{1\}$. Consequently, $H$ is reductive and irreducibility of a pair $(x,y)\in G^2$ implies that $Z_G(H)$ is finite (cf. \cite[Lemma 7.6]{RY00}). Note that, in \cite{BFM02}, irreducible pairs are called rank zero pairs.

\begin{lemma}
Let $G$ be a semisimple algebraic group and let $x,y\in G$ be an irreducible pair satisfying 
\[xy=\omega yx,\quad \omega\in Z(G).\]
Then $x$ and $y$ are semisimple of finite order. Moreover, if $G$ is simply connected then $x$ and $y$ are regular elements.
\end{lemma}
\begin{proof}
Since $Z(G)$ is finite, we may assume that $\omega$ has finite order $m$. Using the relation in the statement, we have
\[x^m =  y\omega^mx^my^{-1} =yx^my^{-1}.  \]
Thus, $x^m\in Z_G(y)$ and similarly, $y^m\in Z_G(x)$. Thus, $x^m$ and $y^m$ are contained in the finite group $ Z_G(H)$, where $H=\overline{\langle x,y\rangle}$, and thus have finite order by \cite[Lemma 7.6]{RY00}. Furthermore, since $x$ and $y$ have finite order they must be semisimple. Moreover, if $G$ is simply connected then $Z_G(x)$ and $Z_G(y)$ are connected \cite[Theorem 2.11]{Hum95} and reductive. Furthermore, $y$ normalises $Z_G(x)$ and thus conjugation by $y$ restricted to $Z_G(x)$ defines an automorphism of $Z_G(x)$ whose fixed-point set is the finite group $Z_G(x)\cap\, Z_G(y) = Z_G(x,y)$.  By \cite[Corollary 10.10]{Ste68}, $Z_G(x)$ is necessarily solvable. Thus, $Z_G(x)$ is a connected reductive solvable group and therefore a torus \cite[Theorem 19.3 (b)]{Hum75}. By symmetry the same argument holds for $Z_G(y)$. Thus, $x$ and $y$ are regular elements.
\end{proof}

The following lemma is the algebraic group version of \cite[Proposition 4.1.1]{BFM02}:

\begin{lemma}\label{lem:conjpairs}
Let $G$ be a semisimple simply connected algebraic group and let
$\omega\in Z(G)$. Suppose that
\[(x_1,y_1),(x_2,y_2)\in G\times G\]
are irreducible pairs satisfying
\[x_i y_i=\omega y_i x_i,\qquad i=1,2.\]
Assume that $x_1$ and $x_2$ are regular semisimple. Then the two pairs
are simultaneously in $G$ if and only if $x_1$ and $x_2$ are
conjugate in $G$.
\end{lemma}
\begin{proof}
The forward direction is trivial. For the converse, we may assume that $x_1=x_2=x\in T$. The goal is to find a $t\in Z_G(x)$ so that $ty_1 t^{-1}=y_2$. Since $y_i$ normalises $T$, the commutator map 
\[\varphi_{y_i}:T\longrightarrow T,\quad t\longmapsto [t,y_i]\] 
is well-defined. Its kernel is
\[\ker\varphi_{y_i} = Z_T(y_i)= Z_G(y_i)\cap T = Z_G(y_i)\cap Z_G(x) = Z_G(y_i,x).\]
Since the pair $(x,y_i)$ is irreducible, $\ker\varphi_{y_i}$ is finite, hence $\dim \im \varphi_{y_i} = \dim T$ and, since the image is closed, we may conclude that $\varphi_{y_i}$ is surjective. Furthermore, one can check that $y_2y_1^{-1}\in Z_G(x)=T$, hence there is a $t\in T$ such that $\varphi_{y_1}(t) = y_2y_1^{-1}$ which is equivalent to $ty_1t^{-1}=y_2$. Since $t\in Z_G(x)$, we have
\[t(x,y_1) t^{-1} = (txt^{-1},ty_1t^{-1})=(x,y_2). \]
\end{proof}

\begin{lemma}
Let $G = \prod_{i=1}^r \SL_{n_i}(\C)$ and let $\omega = \omega_1\cdots \omega_r\in Z(G)$ be a central element where each $\omega_i$ generates $Z(\SL_{n_i}(\C))$. Then all pairs $(x,y)\in G \times G$ satisfying 
\[xy = \omega yx \]
are irreducible.
\end{lemma}
\begin{proof}
Let $U\subset \C^{n_i}$ be a nonzero invariant subspace for the action of $H_i=\overline{\langle x_i,y_i\rangle}$ and $v\in U$ an eigenvector. Since $\omega_i$ has order $n_i$, using the relation $x_iy_i=\omega_i y_i x_i$, it is easy to see that $x_i$ has eigenvectors $\{v, y_iv,\dots, y_i^{n_i-1}v\}$ which are all linearly independent, hence form a basis for $\C^{n_i}$. Thus, $U=\C^{n_i}$ and it follows that each $H_i$ is irreducible and hence so is the pair $(x,y)\in G\times G$.
 
 \end{proof}

\begin{corollary}\label{cor:homogeneous}
Let $p$ be a prime number. Every pair $(x,y)\in \SL_p(\C)\times \SL_p(\C)$ satisfying 
\[xy=\omega yx,\quad \omega\neq 1\in Z(\SL_p(\C))\]
is irreducible and all such pairs are conjugate. In particular, $R_{\omega}(\SL_p(\C))$ is a homogeneous space isomorphic to $\PSL_p(\C)$.
\end{corollary}
\begin{proof}
Let $p$ be a prime number and fix a nontrivial central element $\omega\in Z(\SL_p(\C))$. Let $(x_i,y_i)$ for $i=1,2$ be to pairs in $\SL_p(\C)\times\SL_p(\C)$. Using the defining relation, one sees that the eigenvalues of $x_1$ and $x_2$ are equal, thus $x_1$ and $x_2$ are conjugate and by Lemma~\ref{lem:conjpairs} so are the pairs. If $g\in \Stab_{\SL_p(\C)}((x,y))$, then $g$ is an intertwiner for the restriction to $\overline{\langle x,y\rangle}$ of the standard representation of $\SL_p(\C)$ on $\C^p$. Since $H$ is irreducible, Schur's lemma implies that $g$ is a scalar multiple of the identity, hence $g\in Z(\SL_p(\C))$. Clearly, $Z(\SL_p(\C))\subset \Stab_{\SL_p(\C)}((x,y))$, hence $R_\omega(\SL_p(\C))\simeq \SL_p(\C)/Z(\SL_p(\C))$.
\end{proof}

\begin{theorem}\label{thm:Oka1PSLp}
Let $G=\PSL_p(\C)$, where $p$ is a prime. Then the smooth locus $R(\PSL_p(\C))_{\rm sm}$ of the commuting variety  of $\PSL_p(\C)$ is Oka-$1$.
\end{theorem}
\begin{proof}
Let $q:\SL_p(\C)\to \PSL_p(\C)$ be the quotient map. Then
\[(q\times q)^{-1}(R(\PSL_p(\C))) = \bigsqcup_{\omega\in Z(\SL_p(\C))}R_{\omega}(\SL_p(\C)).\]
Let $\widetilde{G}=\SL_p(\C)$ and $G=\PSL_p(\C)$ and let $\Gamma=\ker(q:\widetilde{G}\to G)$. We have already seen that
\[R(G) = \bigsqcup_{\omega\in \Gamma} R_\omega(G).\]
By Corollary~\ref{cor:homogeneous}, the $R_\omega(\widetilde{G})$, and hence also $R_\omega(G)$, are homogeneous for $\omega\neq 1$. Consequently, the $R_\omega(G)$ with $\omega\neq 1$ are Oka-$1$. The smooth locus $R_{\rm str}(\PSL_p(\C))_{\rm sm}$ is Oka-$1$ by Theorem~\ref{thm:strongOka1}. Therefore all connected components of the full smooth locus are Oka-$1$. The statement of the theorem follows since the Riemann surface in the definition of Oka-$1$ is connected.
\end{proof}

\begin{remark}
The assumption that $p$ is prime in Theorem~\ref{thm:Oka1PSLp} is used to ensure that every non-identity component $R_{\omega}(\PSL_p(\C))$ is a homogeneous space. If  $p$ were not prime, a non-identity element $\omega\in Z(\SL_p(\C))$ need not generate the whole centre. The varieties $R_\omega(\SL_p(\C))$ corresponding to such elements need not be homogeneous spaces. We refer the reader to \cite[\S 4]{BFM02} for the classification of general $\omega$-pairs in the case of compact Lie groups. 
\end{remark}

\section{The Lie algebra case}\label{c:Liealgebra}
In this section we study the holomorphic flexibility properties of the smooth locus of the commuting variety of a semisimple Lie algebra. 

Let $\frakg$ be a rank $r$ complex semisimple Lie algebra with adjoint group $G$. Recall that an element $x\in \frakg$ is \emph{regular} if the dimension of its centraliser $\frakz(x)$ in $\frakg$ is minimal equal to the rank of $\frakg$. The commuting variety of the Lie algebra $\frakg$ is defined as
\[\calC=\calC(\frakg) := \{(x,y)\in \frakg\times \frakg \,|\, [x,y]=0\}.\]

We define two dense open subsets
\[\calC_{i}:=\{(x_1,x_2)\in \calC\,|\,x_i\in \frakg_{\rm reg}\},\quad i=1,\,2.\]
Furthermore, we denote by $\calC_{\rm reg}$ their union
\[\calC_{\rm reg}:=\calC_1\cup \calC_2.\]

\begin{remark}
In \cite[\S 2]{Popov}, the author defines the subset
 \[\calC^\circ = \{(x,y)\in \calC\,|\, \dim (\frakz(x)\cap \frakz(y)) = r\},\]
 also called the regular subset of $\calC$, which is an open subset of the smooth locus $\calC_{\rm sm}$ \cite[Theorem 1.3 (ii)]{Popov}. We remark that the open subset $\calC_{\rm reg} = \calC_1\cup \calC_2$ we defined above is contained in $\calC^\circ$. The distinction already appears in $\sl_4(\C)$ for $x=e_{1,2}+e_{3,4}$ and $y=e_{1,3}+e_{2,4}$ where $e_{i,j}\in \sl_4(\C)$ denotes the standard basis vector whose $i^{\text{th}}$ row and $j^{\text{th}}$ column is $1$ and zero elsewhere. One verifies that $(x,y)\in \calC^\circ$ but neither $x$ nor $y$ is regular elements. Hence, in general there is a proper containment $\calC_{\rm reg}\subset \calC_{\rm sm}$.
\end{remark}
By Chevalley's theorem, we may choose homogeneous generators $p_1, \dots, p_r\in \C[\frakg]^G$. For $x\in \frakg$, the differentials $(dp_i)_x$ may be identified with elements of the dual space $\frakg^*$. Using the non-degeneracy of the Killing form $\kappa$, we have an isomorphism $B:\frakg\to \frakg^*$. Define $F_i(x) \in \frakg$ as the element $F_i(x) := B^{-1}((dp_i)_x)$. Equivalently, $F_i(x)$ is defined by the formula
\[(dp_i)_{x}(v) = \kappa(F_i(x),v).\]

\begin{lemma}\label{lem:Kostantbasis}
For every $x\in \frakg_{\rm reg}$, the $F_1(x),\dots, F_r(x)\in \frakz(x)$ form a basis.
\end{lemma}
\begin{proof}
See \cite[Proposition 1.3]{Kos79} and \cite[Theorem 0.1]{Kos63}.
\end{proof}

\begin{proposition}\label{prop:Liealgbihol}
The map 
\[\Psi:\frakg_{\rm reg}\times \C^r \longrightarrow \calC_1,\quad (x,\lambda_1,\dots, \lambda_r) \longmapsto \left(x,\sum_{i=1}^r\lambda_i F_i(x)\right)\]
is an algebraic isomorphism. 
\end{proposition}

\begin{proof}
The map $\Psi$ is algebraic so we need to show that it has an algebraic inverse. To do so, we construct local inverses on a suitable open cover of $\frakg_{\rm reg}$ and then show they agree on overlaps.

\medskip

Fix a basis $e_1,\dots,e_n$ for $\frakg$ and identify $\frakg\simeq \C^n$. For each $r$-subset $I=\{i_1< \dots < i_r\}$ of the $n$-set $\{1,\dots, n\}$ define the coordinate projection $\pi_I:\C^n \to \C^r$. For each $y\in\frakg$ denote by $y_I$ the image of $y$ under the projection map $\pi_I$ under the above identification. Explicitly, if $y=\sum_{i=1}^n\lambda_ie_i$ then $y_I$ corresponds to
\[\pi_I\left(\sum_{i=1}^n \lambda_i e_i\right) = (\lambda_{i_1},\dots, \lambda_{i_r}).\]
Now, for each $x\in \frakg_{\rm reg}$, let $M(x)$ be the $n\times r$ matrix whose columns are the $F_i(x)$ for $i=1,\dots, r$. That is, 
\[M(x) = [F_1(x)\cdots F_r(x)].\]

Thus, we may write the map $\Psi:\frakg_{\rm reg}\times \C^r\to \calC_1$ as 
\[(x,\lambda_1,\dots,\lambda_r)\longmapsto (x,M(x)(\lambda_1,\dots,\lambda_r)^t).\] 

By Lemma~\ref{lem:Kostantbasis}, for $x\in \frakg_{\rm reg}$ there exists a subset $I=\{i_1< \dots < i_r\}\subset \{1,\dots, n\}$ of $r$ rows of $M(x)$ such that the $r$-by-$r$ sub-matrix 
\[M_I(x):= \pi_I\circ M(x):\C^r\longrightarrow \C^r\]
of $M(x)$ is invertible. For such a choice of $I$, the restriction $\pi_I|_{\frakz(x)}:\frakz(x)\to \C^r$ is an isomorphism. For each $r$-subset $I$ as above, define the open set of $\frakg_{\rm reg}$
\[U_I=\{x\in \frakg_{\rm reg}\,|\, \det(M_I(x))\neq 0\}.\]  
Hence, we have shown that, for each $x\in \frakg_{\rm reg}$ there exists a subset $I$ for which $M_I(x)$ is invertible, and thus the collection of all $\{U_I\}_{I}$ forms an open cover of $\frakg_{\rm reg}$. Corresponding to the open subsets $U_I\subset \frakg_{\rm reg}$, we define
\[\calC_1^{U_I}:=(U_I\times \frakg)\cap \calC =\{(x,y)\in \calC_1\,|\, x\in U_I\}\subset \calC_1.\]
Then, as $I$ varies over all $r$-subsets of $\{1,\dots, n\}$, the open sets $\calC_1^{U_I}$  cover $\calC_1$. For each $\calC_1^{U_I}$ define 
\[\Psi_{U_I}^{-1}:\calC_1^{U_I}\to U_I\times \C^r,\quad (x,y)\mapsto (x,(M_I(x)^{-1}\circ \pi_I)(y)).\]
If $(x,y)\in \calC_1^{U_I\cap U_J}$ for two subsets $I$ and $J$ of $\{1,\dots, n\}$, the restrictions of $\Psi_{U_I}^{-1}$ and $\Psi_{U_J}^{-1}$ to the intersection $\calC_1^{U_I\cap U_J}$ agree. Indeed, since $(x,y)\in \calC_1$, by Lemma~\ref{lem:Kostantbasis} we have  $y\in \frakz(x)=\Span_{\C}\{F_1(x),\dots, F_r(x)\}$. Therefore, there are unique constants $\lambda_1,\dots,\lambda_r$ such that $y=\sum_{i=1}^r\lambda_iF_i(x)$ which under the identification $\frakg\simeq \C^n$ we write as
\[y = M(x)\cdot(\lambda_1,\dots,\lambda_r)^t.\]
Then, in particular, 
\[(M_I(x)^{-1}\circ \pi_I)(y) = (M_I(x)^{-1}\circ \pi_I)(M(x)\cdot(\lambda_1,\dots,\lambda_r)^t) = (\lambda_1,\dots,\lambda_r)^t,\]
and, similarly, $(M_J(x)^{-1}\circ \pi_J)(y) = (\lambda_1,\dots,\lambda_r)^t$. Thus, the $\Psi_{U_{I}}^{-1}$ agree on the overlaps. The $\Psi_{U_I}^{-1}$ define local inverses to $\Psi$ since the restrictions $\Psi_{U_I}:=
\left.\Psi\right|_{U_I\times\mathbb C^r}:U_I\times\mathbb C^r\to \calC_1^{U_I}$
are given by
\[(x,\lambda_1,\dots,\lambda_r)\longmapsto (x,M(x)(\lambda_1,\dots,\lambda_r)^t).\]
Furthermore, since matrix inversion is algebraic, the $\Psi_{U_I}^{-1}$ are algebraic and glue to a global algebraic inverse. 
\end{proof}

\begin{remark}
In \cite[Theorem 1.12]{Popov} the author proves the variety $\calC$ is rational. We give an alternative proof of this fact: The composition $\calC_1\simeq\frakg_{\rm reg}\times \C^r \hookrightarrow  \frakg\times\C^r \simeq \AA^{\dim(\frakg)+r}$ identifies the dense open subset $\calC_1\subset \calC$ with a dense open subset of $\AA^{\dim(\frakg)+r}$, hence $\calC$ is rational.
\end{remark}

\begin{corollary}\label{cor:algellipticity}
The regular commuting variety $\calC_{\rm reg}$ is algebraically elliptic.
\end{corollary}

\begin{proof}
The Zariski open sets $\calC_i$, $i=1,2$ in $\calC_{\rm reg}$ are algebraically elliptic since they are each isomorphic to the elliptic manifold $\frakg_{\rm reg}\times \C^r$. Indeed, by \cite[Theorem 4.12]{Vel72}, $\frakg\setminus\frakg_{\rm reg}$ is a closed subvariety of codimension $3$ in the flexible variety $\frakg\simeq \AA^{\dim(\frakg)}$. By \cite[Theorem 1.1]{FKZ16}, $\frakg_{\rm reg}$ is flexible and hence the product $\frakg_{\rm reg}\times \AA^r$ is algebraically elliptic. Furthermore, algebraic ellipticity is Zariski local, so the union $\calC_{\rm reg}=\calC_1\cup \calC_2$ is algebraically elliptic \cite[Theorem 0.1]{KZ24}.
\end{proof}

\begin{remark}
Applying Corollary~\ref{cor:algellipticity} in the case of $\frakg=\sl_2(\C)$, we see that the smooth locus
\[\calC_{\rm sm} = \calC_{\rm reg} = \calC\setminus \{(0,0)\}\]
is algebraically elliptic (see also \cite[Example 1.10]{Popov}).
\end{remark}

\begin{corollary}
The smooth locus $\calC_{\rm sm}$ is Oka-$1$.
\end{corollary}
\begin{proof}
Set $N= \dim_\C(\calC_{\rm sm})$ and let $E =\calC_{\rm sm}\setminus \calC_1$. Since $\calC_1$ is Zariski dense, $E$ is a proper analytic subvariety with $\dim_\C(E)\leq \dim \calC_{\rm sm}-1$. Consequently, the Hausdorff measure $\calH^{2N-1}_{\calC_{\rm sm},g}(E)$ of $E$ is zero. Since $\calC_1$ is algebraically elliptic, it has the LSAP property and thus $\calC_{\rm sm}$ is Oka-$1$ \cite[Proposition 7.13 (b)]{AF25}.
\end{proof}

\end{document}